\documentclass[aap,MSNbibl,dvips]{arximspdf}
\usepackage{graphicx}
\doi{10.1214/11-AAP820} 
\volume{23}
\issue{1}
\pubyear{2013}
\firstpage{99}
\lastpage{144}

\makeatletter

\newtheorem{theorem}{Theorem}[section]

\newproclaim{dfn}[theorem]{Definition}
\newproclaim{rem}{Remark}

\newtheorem{prop}[theorem]{Proposition}

\newcommand{\cal}{\mathcal}

\newcommand{\lbd}{\lambda}
\newcommand{\eps}{\epsilon}
\newcommand{\vareps}{\varepsilon}
\newcommand{\tendinfty}{\rightarrow\infty}

\newcommand{\supp}{\mathfrak{s}}
\newcommand{\ZZ}{\mathbb{Z}}
\newcommand{\NN}{\mathbb{N}}
\newcommand{\RR}{\mathbb{R}}
\newcommand{\PP}{\mathbb{P}}
\newcommand{\EE}{\mathbb{E}}

\newcommand{\indic}[1]{{\mathbf{1}}_{\{#1\}}}
\newcommand{\indicbis}[1]{{\mathbf{1}}_{#1}}

\newcommand{\kil}{\mathrm{k}}

\makeatother

\begin{document}
\begin{frontmatter}

\title{The coalescent point process of branching trees}
\runtitle{The coalescent point process of branching trees}

\begin{aug}
\author[A]{\fnms{Amaury} \snm{Lambert}\thanksref{t1}}
\and
\author[B]{\fnms{Lea} \snm{Popovic}\corref{}\thanksref{t2}\ead[label=e1]{lpopovic@mathstat.concordia.ca}}
\runauthor{A. Lambert and L. Popovic}
\affiliation{UPMC Univ Paris 06 and Concordia University}
\address[A]{\'{E}quipe Probabilit\'{e}s, Statistiques\\
\quad \& Biologie\\
Laboratoire de Probabilit\'{e}s\\
\quad \& Mod\`{e}les Al\'{e}atoires\\
UPMC Univ Paris 06\\
Case courrier 188\\
4, Place Jussieu\\
75252 Paris Cedex 05\\
France} 

\address[B]{Department of Mathematics and Statistics\\
Concordia University\\
1455 De Maisonneuve Blvd West\\
Montreal, H3G IM8\\
Canada}
\end{aug}

\thankstext{t1}{Supported by the Projects MAEV 06-BLAN-3 146282 and MANEGE
09-BLAN-0215 of ANR (French national research agency).}

\thankstext{t2}{Supported by the University Faculty award and Discovery
grant of NSERC (Natural
sciences and engineering council of Canada).}

\received{\smonth{1} \syear{2011}}
\revised{\smonth{10} \syear{2011}}

%
\begin{abstract}
We define a doubly infinite,
monotone labeling of Bienaym\'{e}--Galton--Watson (BGW) genealogies.
The genealogy of the current generation backwards in time is uniquely
determined by the coalescent point process $(A_i$; $i\ge1)$, where $A_i$
is the coalescence time between individuals $i$ and $i+1$. There is a
Markov process of point measures $(B_i;i\ge1)$ keeping track of more
ancestral relationships, such that $A_i$ is also the first point mass
of $B_i$.

This process of point measures is also closely related to an
inhomogeneous spine decomposition of the lineage of the first surviving
particle in generation $h$ in a planar BGW tree conditioned to survive
$h$ generations. The decomposition involves a point measure $\rho$
storing the number of subtrees on the right-hand side of the spine.
Under appropriate conditions, we prove convergence of this point
measure to a point measure on $\RR_+$ associated with the limiting
continuous-state branching (CSB) process. We prove the associated
invariance principle for the coalescent point process, after we
discretize the limiting CSB population by considering only points with
coalescence times greater than $\vareps$.

The limiting coalescent point process $(B^\vareps_i;i\ge1)$ is the
sequence of depths greater than $\vareps$ of the excursions of the
height process below some fixed level. In the diffusion case, there are
no multiple ancestries and (it is known that) the coalescent point
process is a Poisson point process with an explicit intensity measure.
We prove that in the general case the coalescent process with
multiplicities $(B^\vareps_i;i\ge1)$ is a Markov chain of point masses
and we give an explicit formula for its transition function.

The paper ends with two applications in the discrete case. Our results
show that the sequence of $A_i$'s are i.i.d. when the offspring
distribution is linear fractional. Also, the law of Yaglom's
quasi-stationary population size for subcritical BGW processes is
disintegrated with respect to the time to most recent common ancestor
of the whole population.
\end{abstract}

%

\end{frontmatter}

\section{Introduction}
The idea of describing the backward genealogy of a population is
ubiquitous in population genetics. The most popular piece of work on
the subject is certainly~\cite{Ki}, where the \textit{standard
coalescent} is introduced, and shown to describe the genealogy of a
finite sample from a population with large but constant population
size. Coalescent processes for branching processes cannot be
characterized in the same way, since, for example, they are not
generally Markov as time goes backwards, although for stable
continuous-state branching (CSB) processes the genealogy can be seen as
the time-change of a Markovian coalescent~\cite{7}.

The present paper relies on the initial works~\cite{AP,P}, which
focused on the coalescent point process for critical birth--death
processes and the limiting Feller diffusion. They have been extended to
noncritical birth--death processes~\cite{G} and more generally to
homogeneous, binary Crump--Mode--Jagers processes~\cite{L2}. In all
these references, simultaneous births were not allowed, since then the
genealogical process would have to keep memory of the multiplicity of
all common offspring of an ancestor. The problem of sampling a
branching population or a coalescent point process has received some
attention~\cite{L0,L3,S}, but no consistent sampling in the standard
coalescent has been given so far (except Bernoulli sampling of leaves).
Our goal was to define a coalescent point process for arbitrary
branching processes, that is both simple to describe in terms of its
law, and allows for finite sampling in a consistent way. That is, the
coalescent process for samples of size $n\ge1$ are all embedded in the
same object.

A different way of characterizing the genealogy of a branching
population alive at some fixed time is with a \textit{reduced tree},
first studied in~\cite{FSS} and generalized in~\cite{DLG}, Section 2.7,
to CSB processes. To construct it, a starting date $0$ and a finishing
date $T$ are specified, and a reduced branching tree started at $0$ and
conditioned to be alive at $T$ is defined by erasing the points without
alive descendants at time $T$. Instead of directly displaying the
coalescence times as a sequence running over the current population
size (as is our goal), this approach characterizes the transition
probabilities of the reduced branching process by tracking the
population size in time with an inhomogeneous Markov process on $[0,T]$
taking values in the set of integers for all times in $[0,T)$.
Unfortunately, this construction does not allow for a consistent way of
sampling the individuals alive at time $T$.

We use a different approach and construct a coalescent process with
multiplicities for the genealogy of some random population, when the
forward time genealogy is produced by a general branching process,
either discrete or continuous-state. Our main goal is to give a simple
representation for this process, and describe its law in a manner that
would be easy to use in applications.

\begin{figure}

\includegraphics{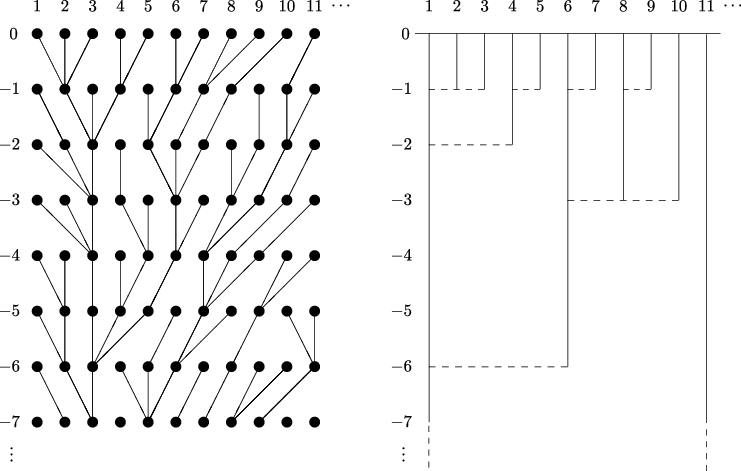}

\caption{A doubly infinite embedding of quasi-stationary Bienaym\'
{e}--Galton--Watson genealogies. Horizontal labels are individuals and
vertical labels are generations. Left panel: the complete embedding of
an infinite BGW tree. Right panel: coalescence times between
consecutive individuals of generation 0. Here, the coalescent point
process takes the value $(1, 1, 2, 1, 6, 1, 3, 1, 3, \ldots)$.}
\label{figcoalBGW}
\end{figure}

We are interested in having an arbitrarily large population at the
present time, that we think of as generation 0, arising from a general
branching process, originating at an unspecified arbitrarily large time
in the past.
In order to keep track of genealogies of individuals in the present
population, for discrete Bienaym\'{e}--Galton--Watson (BGW) branching
processes we use a representation that labels individuals at each
generation in such a way that lines of descent of the population at any
time do not intersect (see Figure~\ref{figcoalBGW}). This leads to a
monotone planar embedding of BGW trees such that each individual is
represented by $(n,i)$ where $n\in\ZZ$ is the generation number and
$i\in\NN$ is the individual's index in the generation. We consider BGW
trees which are doubly infinite, allowing the number of individuals
alive at present time to be arbitrarily large, and considering an
infinite number of generations for their ancestry back in time. This
monotone representation can also be extended to the case of
continuous-state branching processes (later called CSB processes), in
terms of a linearly ordered planar embedding of the associated
genealogical $\RR$-trees. This is in some way implicit in many recent
works: our embedding is the exact discrete analogue of the flow of
subordinators defined in~\cite{BLG}; also, in the genealogy of CSB
processes defined in~\cite{LGLJ}, the labeling of individuals at a
fixed generation is analogous to the successive visits of the height
process at a given level.

From our discrete embedding of the population, a natural consequence is
that \textit{coalescence times}, or times backwards to the most recent
common ancestor, between individuals $k$ and $\ell$ in generation $0$,
are represented by the maximum of the values $A_i, i=k,\ldots, \ell-1$,
where $A_i$ is the coalescence time between individuals $i$ and $i+1$.
The genealogy back in time of the present population is then uniquely
determined by the process $(A_i; i\ge1)$, which we call the
\textit{branch lengths} of the coalescent point process. In the
continuous-state branching process setting, the $A_i$'s are given by
depths of the excursions of the height process below the level
representing the height of individuals in the present population.

In general, the process $(A_i; i\ge1)$ is not Markovian and its law is
difficult to characterize. Our key strategy is to keep track of
multiple ancestries to get a Markov process, thus constructing a
\textit{coalescent point process with multiplicities} $(B_i;i\ge1)$.
Intuitively, for each $i\ge1$, $B_i$ encodes the relationship of the
individual $i+1$ to the infinite spine of the first present individual,
by recording the nested sequence of subtrees that form that ancestral
lineage linking $i+1$ to the spine. The value of $B_i$ is a point mass
measure, where each point mass encodes one of these nested subtrees:
the level of the point mass records the level back into the past at
which this subtree originated, while the multiplicity of the point mass
records the number of subtrees with descendants in individuals $\{i+1,
i+2,\ldots\}$ emanating at that level which are embedded on the
right-hand side of the ancestral link of $i+1$ to the
spine.\looseness=-1

Formally, for both BGW and the CSB population, we proceed as follows.
We define a process with values in the integer-valued measures on $\NN$
or $\RR_+$ (BGW or CSB, resp.) in a recursive manner, in such a way
that $(B_j;1\le j\le i)$ will give a complete record of the ancestral
relationship of individuals $1\le j\le i$. Start with the left-most
individual in the present population, $i=1$, and let $B_1$ have a point
mass at $n$, where $-n$ is the generation of the last common ancestor
of individuals $1$ and $2$, and have multiplicity $B_1(\{n\})$ equal to
the number of times this ancestor will also appear as the last common
ancestor for individuals ahead, with index $i\ge1$. We then proceed
recursively, so that at individual $i$ we take the point masses in
$B_{i-1}$; we first make an adjustment in order to reflect the change
in multiplicities due to the fact that the lineage of individual $i$ is
no longer considered to be an individual ahead of the current one. We
then add a record for the last common ancestor of individuals $i$ and
$i+1$. If this is an ancestor that has already been recorded in $B_
{i-1}$, we just let $B_i$ be the updated value of $B_{i-1}$. Otherwise,
we let $B_i$ be the updated $B_{i-1}$ plus a new point mass at $n$, if
$-n$ is the generation of this new last common ancestor, whose\vadjust{\goodbreak}
multiplicity $B_i(\{n\})$ is the number of times this new ancestor will
appear as the last common ancestor for individuals ahead, with index
$\ge i+1$ (e.g., in Figure~\ref{figcoalBGW} we will have $B_1=2\delta
_{1}, B_2=\delta_{1}, B_3=\delta_{2}, B_4=\delta_{1}, B_5=\delta_{6},
B_6=\delta_{1}, B_7=2\delta_{3}, B_8=\delta_{1}+\delta_{3},
B_9=\delta_{3}$).

Once we have constructed the coalescent point process with
multiplicities (in either the BGW or CSB case), we will show that $A_i$
can be recovered as the location of the nonzero point mass in $B_i$
with the smallest level, that is, $A_i=\inf\{n\dvtx  B_i(\{n\})\neq0\}$
(e.g., in Figure~\ref{figcoalBGW} we have $A_1=1, A_2=1, A_3=2$, \mbox{$A_4=1$},
$A_5=6, A_6=1, A_7=3, A_8=1, A_9=3$). More importantly, we prove that
$(B_i;i\ge1)$ is a Markov process, and that, when going from
individual $i$ to $i+1$, the transitions decrease the multiplicity of
the point mass of $B_i$ at level $\inf\{n\dvtx  B_i(\{n\})\neq0\}$ by 1,
and that, with a specified probability, a new random point mass is
added at a random level that must be smaller than the smallest nonzero
level in the updated version of $B_i$.

In order to use this construction for both BGW and the CSB population,
our take on what constitutes the sequence of present individuals has to
be different for the continuous CSB population from the simple one for
the discrete BGW population. Since in the CSB case the present
population size is not discrete, and there is an accumulation of
immediate ancestors at times arbitrarily close to the present time, we
have to discretize the present population by considering only the
individuals whose last common ancestor occurs at a time at least an
$\vareps$ amount below the present time, for an arbitrary $\vareps>0$.
We will later show that we can obtain the law of the coalescent point
process with multiplicities for the CSB population as a limit of a
sequence of appropriately rescaled coalescent point processes with
multiplicities for the BGW population for which we also use the same
discretization process of the present population. 
At first it may seem surprising that these coalescents with
multiplicities are Markov processes over the set of all or the $\vareps
$-discretized set (BGW and CSB case, resp.) individuals at present
time. Below we intuitively explain why this is the case by describing
the two approaches for constructing them.

In the discrete case, we start by giving a related, easier to define
process $(D_i$; \mbox{$i\ge1)$}, taking values in the integer-valued sequences,
whose first nonzero term is also at level $A_i$. For each $i$, the
sequence $(D_i(n),n \ge1)$ gives the number of younger offshoots at
generation $-n$ embedded on the right-hand side of the ancestor of $i$.
The trees sprouting from the younger offshoots are independent, and the
law of a tree sprouting from a younger offshoot at generation $-n$ has
the law of the BGW tree conditioned to survive at least $n-1$
generations. It turns out that $(D_i;i\ge1)$ is Markov, and we
construct $(B_i;i\ge1)$ from it, show that it is Markov as well and
give its transition law.
In order to be able to deal with the conditioning of younger offshoots
in a way that allows us later to pass to the limit, we introduce an
integer-valued measure $\rho$ that takes the ancestor of individual $1$
at generation $-n$ and records as $\rho(\{n\})=\rho_n$ the number of
all of its younger offshoots embedded on the\vadjust{\goodbreak} right-hand side, and we
call it the \textit{great-aunt measure}. This measure gives a spine
decomposition of the first survivor (individual $i=1$) in such a way
that the law of the trees sprouting from the younger offshoots are
still independent, but are no longer conditioned on survival.

In the continuous-state case, the great-aunt measure will be a measure
$\rho^0$ on~$\RR_+$, and we will be able to define the genealogy thanks
to independent CSB processes starting from the masses of the atoms of
$\rho^0$. In the subcritical and critical cases, this can be done using
a single path of a L\'{e}vy process with no negative jumps and Laplace
exponent $\psi$. In the supercritical case a concatenation of excursion
paths will have to be used. Using the continuous great-aunt measure
$\rho^0$, we characterize the genealogy of an infinite CSB tree with
branching mechanism $\psi$ via the height function $H^\star$ whose
value at an individual in the population can be decomposed into the
level on the infinite spine at which the subtree containing this
individual branches off, and the relative height of this individual
within its subtree. We construct $(B^\vareps_i;i\ge1)$ from the height
process $H^\star$.
Discretizing the population by considering only the points whose
coalescence times are greater than some fixed $\varepsilon>0$,
translates into considering only the excursions of $H^\star$ from level
$0$ with a depth greater than $-\vareps$. From these excursions we can
obtain the process $((A^\vareps_i, N^\vareps_i);i\ge1)$ where
$A^\vareps_i$ is the depth of the $i$th such excursion and $N^\vareps
_i$ is the number of future excursions with the exact same depth. It
turns out that a specific functional of this process is Markov, and we
construct $(B^\vareps_i;i\ge1)$ from it, show that it is Markov as
well and give its transition law.

\begin{figure}

\includegraphics{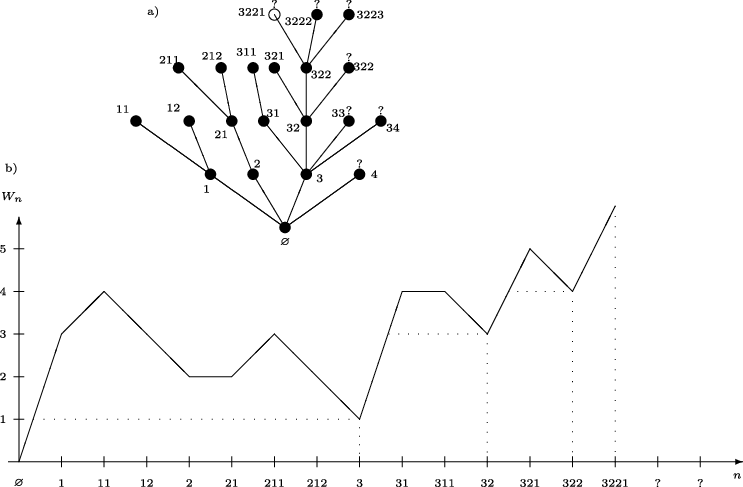}

\caption{\textup{(a)} A Bienaym\'{e}--Galton--Watson tree conditioned to have
alive individuals at generation~4. The empty circle is the first such
individual (3221) in the lexicographical order. Descendants of
individuals with greater rank are only indicated by a question mark;
\textup{(b)} the associated random walk $W$ killed after the visit $\sigma_4$ of
the first individual $x_4=3221$ with height 4. Records of the future
infimum at $\sigma_4$ are shown by dotted lines. The record times are
exactly the visits to ancestors $322, 32, 3$ and $\varnothing$ of 3221.
The pairs of overshoots and undershoots of $W$ across the future
infimum at those times (in this order) are $({\alpha}_1, {\rho
}_1)=(0,2)$, $({\alpha}_2, {\rho}_2)=(1,1)$, $({\alpha}_3, {\rho
}_3)=(1,2)$ and $({\alpha}_4, {\rho}_4)=(2,1)$. }
\label{figrandomwalk}
\end{figure}

Now, in order to prove convergence in law of the discretized version of
$(B_i$; $i\ge1)$ to $(B^\vareps_i;i\ge1)$, we make use of the great-aunt
measures $\rho$ and $\rho^0$ describing the spine decomposition of the
first surviving individual.
In the discrete case, the spine decomposition $(\rho(\{n\}); 0\le n\le
h)$ truncated at level $h$ has the same law as the spine decomposition
of a planar embedding of a BGW tree conditioned on surviving up to
generation $h$, where the first survivor is defined in the usual
depth-first search order (see Figure~\ref{figrandomwalk}). Using a
well-known random walk representation of BGW trees~\cite{BK}, we then
prove (under usual conditions ensuring the convergence of the random
walk to a spectrally positive L\'{e}vy process $X$ with Laplace
exponent~$\psi$) the convergence (in the vague topology) of the
great-aunt measure $\rho$ to the measure $\rho^0$ on $\RR_+$ defined by
\[
\rho^0(dx):= \beta \,dx+ \sum_{t\dvtx \Delta_t>0}
\Delta_t\delta_t(dx),
\]
where $\beta$ is the Gaussian coefficient of $X$, and $(t,\Delta_t)$ is
a Poisson point measure with intensity
\[
\biggl(e^{v(t)r} \int_{(r,\infty)} e^{-v(t)z}
\Lambda(dz) \biggr) \,dt \,dr,\qquad r>0,
\]
where $\Lambda$ is the L\'{e}vy measure of $X$, and $v$ the inverse of
the decreasing function\vadjust{\goodbreak} $\lbd\mapsto\int_{[\lbd,+\infty)}du/\psi(u)$.
This is exactly the same measure $\rho^0$ we obtain from the spine
decomposition in the continuous case. Since the transition law of
$(B_i;i\ge1)$ can be represented as a functional of $\rho$ and the
transition law of $(B^\vareps_i;i\ge1)$ as a functional of $\rho^0$,
this will lead to our claim.

Finally, in the very last section, we give two simple applications of
our results in the discrete case. First, we prove that in the
linear-fractional case, the coalescent point process is a sequence of
i.i.d. random variables. Related results can be found in~\cite{Ran}.
Second, in the subcritical case, we use the monotone embedding to
display the law, in quasi-stationary state, of the total population
size (Yaglom's limit) jointly with the time to most recent common
ancestor of this whole population.

\section{Doubly infinite embedding and the coalescent point
process}\label{sec2}

\subsection{The discrete model}

We will start off with a monotone planar embedding of an infinitely old
BGW process with arbitrarily large population size.\vadjust{\goodbreak} Let $(\xi(n,i))$ be
a doubly indexed sequence of integers, where $n$ is indexed by $\ZZ=\{
\ldots, 0,-1,-2,-3,\ldots\}$ and $i$ is indexed by $\NN=\{1,2,3,\ldots
\}$.
The index pair $(n,i)$ represents the $i$th individual in generation
$n$, and $\xi(n,i)$ provides the number of offspring of this individual.

We endow the populations with the following genealogy.
Individual $(n,i)$ has mother $(n-1,j)$ if
\[
\sum_{k=1}^{j-1} \xi(n-1,k) < i \le\sum
_{k=1}^{j} \xi(n-1,k).
\]
From now on, we focus on the ancestry of the population at time 0, that
we will call \textit{standing population}, and we let $\mathfrak{a}_i(n)$
denote the index of the ancestor of individual $(0,i)$ in generation
$-n$. In particular, $\mathfrak{a}_i(1):=\min\{j\ge1\dvtx  \sum_{k=1}^{j}
\xi(-1,k)\ge i\}$.
Our main goal is to describe the law of the times of coalescence
$C_{i,j}$ of individuals $(0,i)$ and $(0,j)$, that is,
\[
C_{i,j}:=\min\bigl\{n\ge1\dvtx  \mathfrak{a}_i(n)=
\mathfrak{a}_j(n)\bigr\},
\]
where it is understood that $\min\varnothing= +\infty$. Defining
\[
A_i:=C_{i,i+1},
\]
it is easily seen that by construction, for any $i\le j$, $C_{i,j}=
\max
\{A_i,A_{i+1},\ldots,\break A_{j-1}\}$. Thus, the sequence $A_1, A_2,\ldots$
contains all the information about the (unlabeled) genealogy of the
current population and is called \textit{coalescent point process}, as in
\cite{P} (see Figure~\ref{figcoalBGW}).

We assume that there is a random variable $\xi$ with values in $\ZZ_+=\{
0,1,2,\ldots\}$ and probability generating function (p.g.f.) $f$, such
that all random variables (r.v.s) $\xi(n,i)$ are i.i.d. and
distributed as $\xi$. As a consequence, if $Z^{(n,i)}(k)$ denotes the
number of descendants of $(n,i)$ at generation $n+k$, then the
processes $(Z^{(n,i)}(k);k\ge0)$ are identically distributed Bienaym\'
{e}--Galton--Watson (BGW) processes starting from 1. We will refer to
$Z$ as some generic BGW process with this distribution and we will
denote by $f_n$ the $n$th iterate of $f$ and by $p_{n}=1-f_{n}(0)$ the
probability that $Z_n\not=0$.
We need some extra notation in this setting.

We define $\zeta_n$ as the number of individuals at generation 1 having
alive descendants at generation $n$.
In particular, $\zeta_n$ has the law of
\[
\sum_{k=1}^\xi\eps_k,
\]
where the $\eps_k$ are i.i.d. Bernoulli r.v.s with success probability
$p_{n-1}$, independent of $\xi$. We also define $\zeta_n'$ as the r.v.
distributed as $\zeta_n-1$ conditional on $\zeta_n\not=0$.

\subsection{The coalescent point process: Main results}\label{subsecdiscrcoal}
Recall that $A_i$ denotes the time of coalescence of individuals
$(0,i)$ and $(0,i+1)$, or $i$th \textit{branch length of the coalescent}.
To describe the law of $(A_i;i\ge1)$ we need an additional process.\vadjust{\goodbreak}

Let $D_i(n)$ be the number of daughters of $(-n, \mathfrak{a}_i(n))$,
distinct from $(-n+1$, $\mathfrak{a}_i(n-1))$, having descendants in $\{
(0,j);j\ge i\}$. In other words, $D_i(n)$ is the number of younger
surviving offshoots of $(-n,\mathfrak{a}_i(n))$ not counting the
lineage of $(0,i)$ itself. Letting
\[
{\cal D}_i({n}):=\bigl\{\mbox{daughters of $\bigl(-n,
\mathfrak{a}_i(n)\bigr)$ with descendants in $\bigl\{(0,j);j\ge i
\bigr\}$}\bigr\},
\]
we have
\[
D_i(n)=\#{\cal D}_i({n})-1.
\]
We set $D_0(n):=0$ for all $n$. We now provide the law of this process
and its relationship to $(A_i;i\ge1)$. We also set $A_0:=+\infty$,
which is in agreement with $\min\varnothing=+\infty$.
%
\begin{theorem}
\label{thmcpp}
The $i$th branch length is a simple functional of $(D_i;i\ge0)$
\[
A_i=\min\bigl\{n\ge1\dvtx  D_i(n)\not=0\bigr\}.
\]
In addition, the sequence-valued chain $(D_i;i\ge0)$ is a Markov chain
started at the null sequence. For any sequence of nonnegative integers
$(d_n)_{n\ge0}$, the law of $D_{i+1}$ conditionally given $D_i(n)=d_n$
for all $n$, is given by the following transition. We have
$D_{i+1}(n)=d_n$ for all $n>A_i$, $D_{i+1}(A_i)=d_{A_i}-1$ and the
r.v.s $D_{i+1}(n)$, $1\le n<A_i$, are independent, distributed as
$\zeta_n'$.
In particular, 
the law of $A_1$ is given by
\begin{eqnarray*}
\PP(A_1> n)&=& \prod_{k=1}^n\PP\bigl(
\zeta_{k}'=0\bigr) = \frac{1}{1-f_{n}(0)}
\prod_{k=1}^nf' \bigl(f_{k-1}(0)
\bigr) \\
&=& \frac{f_n'(0)}{1-f_{n}(0)}= \PP(Z_n=1\mid Z_n\not=0).
\end{eqnarray*}
\end{theorem}
\begin{pf}
The following series of equivalences proves that $A_i$ is the level of
the first nonzero term of the sequence $D_i$.
\begin{eqnarray*}
A_i>n \quad&\Leftrightarrow&\quad \forall k\le n, \mathfrak{a}_i(k)
\not=\mathfrak{a}_{i+1}(k)
\\
&\Leftrightarrow&\quad \forall k\le n, \forall j>i, \mathfrak{a}_i(k)
\not=\mathfrak{a}_j(k)
\\
&\Leftrightarrow&\quad \forall k\le n, \bigl(-k, \mathfrak{a}_i(k)\bigr)
\mbox{ has no descendants in }\bigl\{(0,j)\dvtx j>i\bigr\}
\\
&\Leftrightarrow&\quad \forall k\le n, D_i(k)=0.
\end{eqnarray*}
Thanks to this last result, we get
\[
\mathfrak{a}_i(A_i)=\mathfrak{a}_{i+1}(A_i)
\quad\mbox{and}\quad\mathfrak{a}_i(A_i-1)\not=
\mathfrak{a}_{i+1}(A_i-1).
\]
In particular, $(-A_i+1,\mathfrak{a}_i(A_i-1))$ has no descendants in
$\{(0,j)\dvtx j\ge i+1\}$, so that
\[
{\cal D}_{i+1} (A_i) = {\cal D}_i
(A_i) \setminus\bigl\{ \bigl(-A_i+1, \mathfrak
{a}_i(A_i-1)\bigr)\bigr\}
\]
and $\#{\cal D}_{i+1} (A_i)=\#{\cal D}_i (A_i)-1$, that is,
$D_{i+1}(A_i)=D_i(A_i)-1$.\vadjust{\goodbreak}


Let us deal with the case $n>A_i$. By definition of $A_i$, $\mathfrak
{a}_i(n)=\mathfrak{a}_{i+1}(n)$ for any $n\ge A_i$. As a consequence,
for any $n>A_i$, each daughter of $(-n, \mathfrak{a}_{i+1}(n))$ has
descendants in $\{(0,j);j\ge i+1\}$ iff she has descendants in $\{
(0,j);j\ge i\}$. In other words, ${\cal D}_{i+1}(n)={\cal D}_i(n)$ and
$D_{i+1}(n)=D_i(n)$.


Now we deal with the case $n<A_i$. Set
\[
E_i:=\bigl\{D_j(n)=d_{j,n}, \forall n\ge1, j
\le i\bigr\},
\]
where the $(d_{j,n})_{n\ge1, j\le i}$ are fixed integer numbers, and
let ${A_i}$ be the value of the coalescence time of $i$ and $i+1$
conditional on $E_i$, that is, ${A_i}=\min\{n\ge1\dvtx d_{i,n}\not=0\}$.
Now let ${\cal T}(n,i)$ be the tree descending from $(-n,\mathfrak
{a}_i(n))$ and set
\[
I(n,i):=\min\bigl\{j\le i\dvtx  \mathfrak{a}_j(n)=\mathfrak{a}_i(n)
\bigr\}.
\]
Observe that the (unlabeled) tree ${\cal T}(n,i)$ has the law of a BGW
tree conditioned on having at least one descendant at generation 0. 
Now, because the event $E_i$ only concerns the descendants of daughters
of ancestors of $(-{A_i}+1, \mathfrak{a}_{i+1}({A_i}-1))$, the law of
${\cal T}({A_i}-1,i+1)$ conditional on $E_i$ is still the law of a BGW
tree conditioned on having at least one descendant at generation 0.

Also notice that conditional on $E_i$, $I({A_i}-1,i+1)=i+1$. As a
consequence, $I(n,i+1)=i+1$ for any $n<{A_i}$, so that for any
$n<{A_i}$, conditional on $E_i$,
\[
{\cal D}_{i+1}(n):=\bigl\{\mbox{daughters of $\bigl(-n,
\mathfrak{a}_{i+1}(n)\bigr)$ with descendants at generation 0}\bigr\}.
\]
The result follows recalling the law of ${\cal T}({A_i}-1,i+1)$
conditional on $E_i$.
\end{pf}

Recall that $D_0$ is a null sequence, $A_0=+\infty$, and that
$D_1(n)\stackrel{d}{=}\zeta'(n)$ for $n\ge1$. This infinite sequence
of values in $\{D_1(n),n\ge1\}$ contains information on the ancestral
relationship of the individual $(0,1)$ and an arbitrarily large number
of individuals in the standing population going back arbitrarily far
into the past.
Using the process $(D_i;i\ge1)$ we next define a sequence $(B_i;i\ge
1)$ of finite point mass measures which contains the minimal amount of
information needed to reconstruct $A_1, A_2,\ldots$ while remaining
Markov. We do this for two reasons. First, if we are only interested in
the ancestral relationship of finitely many individuals in the standing
population, there is no need to keep track of an infinite sequence of
values. Second, when we consider a rescaled limit of the BGW process to
a CSB process, we will need to work with a sparser representation.

The main distinction between these two processes is that while $D_i$
contains information on the ancestral relationship of $(0,i)$ and
$(0,j)$ for both $1\le j\le i+1$ and $j> i+1$, $B_i$ will only contain
information about the ancestral relationship of $(0,i)$ and $(0,j)$ for
$1\le j\le i+1$. In other words, if, say, $\max\{A_1,\ldots, A_i\}=n$,
then $B_i(m)=0$ for all $m>n$. Moreover, $B_i$ is defined directly from
$D_i$ by letting $B_i(\{m\})=D_i(m)$ for all $m\le n$. In particular,
$B_1$ will have a single nonzero point mass at level $A_1$, and $B_1(\{
A_1\})=D_1(A_1)$.

We are now ready to define $(B_i;i\ge1)$ which we call the
\textit{coalescent point process with multiplicities}. Let
$(B_i;i\ge0)$ be a sequence of finite point measures, started at the
null measure, defined from $(D_i;i\ge0)$ recursively as follows. For
any point measure $b=\sum_{n\ge1}b_{n}\delta_n$, let $\supp(b)$ denote
the minimum of the support of~$b$, that is,
\[
\supp(b):=\min\{n\ge1\dvtx  b_{n}\neq0\}
\]
with the convention that $\supp(b)=+\infty$ if $b$ is the null measure.
If $B_i=\break\sum_{n\ge1}b_{n,i}\delta_n$ for some $\{b_{n,i}\}_{n\ge
1}\in
\mathbb{N}$, let
\[
a_{1,i}:= \supp(B_i),\qquad B_i^*:=
B_i-\delta_{a_{1,i}},\qquad a_{1,i}^*:= \supp
\bigl(B_i^*\bigr).
\]
%
Then, define
\[
B_{i+1}:= \cases{ B_i^* + D_{i+1}(A_{i+1})
\delta_{A_{i+1}}, &\quad if $A_{i+1}<a_{1,i}^*$ and
$A_{i+1}\neq a_{1,i}$,
\vspace*{1pt}\cr
B_i^*, &\quad otherwise. }
\]
Note that by Theorem~\ref{thmcpp} we have $D_{i+1}(A_i)=D_i(A_i)-1$,
so by this definition $B_{i+1}=\sum_{n\ge1}b_{n,i+1}\delta_n$, where
$\{b_{n,i+1}\}_{n\ge1}\in\mathbb{N}$ satisfies
\[
b_{n,i+1}:= \cases{ D_{i+1}(n), &\quad if $n$ is such that
$b_{n,i}\neq0$,
\cr
D_{i+1}(A_{i+1}), &\quad if
$n=A_{i+1}$ and $A_{i+1}<a_{1,i}^*$ and
$A_{i+1}\neq a_{1,i}$,
\cr
0, &\quad for all other $n$. }
\]
%

Roughly, $(B_i;i\ge1)$ records the ancestral information in the
following way. $B_0$ is a null measure, and $B_1$ will contain a single
point mass $B_1=D_1(A_1)\delta_{A_1}$ recording the coalescence time
$A_1$ for individuals $(0,1)$ and $(0,2)$ in the location of its point
mass. Since the last common ancestor of $(0,1)$ and $(0,2)$ may also be
the last common ancestor of $(0,1)$, $(0,2)$ and some other individual
$(0,j)$ for $j\ge2$, the multiplicity of this point mass will record
the number of its future appearances in the coalescent point process.
Recursively in every step, say $i+1$, this record of point masses will
need to be updated from $B_i$ to $B^*_i$, since by moving one
individual to the right we are no longer recording the last common
ancestor of two previous individuals, and the number of future
appearances of their last common ancestor by definition goes down by
$1$. In addition, at step $i+1$ we also need to record the coalescence
time for the last common ancestor of $(0,i+1)$ and $(0,i+2)$ with the
number of its future appearances. This is done by taking the updated
version $B^*_i$ and adding a new point mass $D_{i+1}(A_{i+1})\delta
_{A_{i+1}}$ to create $B_{i+1}$. Because of the monotone embedding of
the BGW as a planar tree, it is not possible for the level $A_{i+1}$ of
this new point mass to be greater than any of the common ancestors of
$(0,j)$ and $(0,k)$ for $1\le j<k\le i+1$ unless the multiplicities of
these ancestors are depleted at this step and no longer appear in the
updated $B^*_i$ (so $A_{i+1}\le a^*_{1,i}$) as will be seen in the
proof of the next theorem. Moreover, if the last common ancestor of
$(0,i+1)$ and $(0,i+2)$ is the last common ancestor of $(0,i+1)$ and
$(0,k)$ for some $1\le k\le i+1$, then the common ancestry of $(0,i+2)$
and $(0,i+1)$ will be counted in the multiplicity of the mass at
$A_{i+1}$ in $B_i$ and the updated version $B^*_i$ will have nonzero
multiplicity at $A_{i+1}$ (so $A_{i+1}=a^*_{1,i}$) and there will be no
need to add a new point mass to create $B_{i+1}$. In addition, if this
ancestor is also the last common ancestor of $(0,i)$ and $(0,i+1)$ (so
$A_{i+1}=a_{1,i}$), the count for this ancestor cannot be 1 in $B_i$,
so in this case we have $A_{i+1}= a_{1,i}=a^*_{1,i}$ and there is no
need to add a new point mass in $B_{i+1}$.

We now provide the law of this point measure process and its
relationship to $(A_i;i\ge1)$.
%
\begin{theorem}
\label{thmcpmp}
The $i$th branch length is the smallest point mass in $B_i$
\[
A_i=\supp(B_i)=\min\bigl\{n\ge1\dvtx  B_i
\bigl(\{n\}\bigr)\neq0\bigr\}.
\]
In addition, the sequence of finite point measures $(B_i; i\ge0)$ is a
Markov chain started at the null measure, such that for any finite
point measure $b=\sum_{n\ge1}b_n\delta_{n}$, with $b_n\in\mathbb
{N}\cup\{0\}$, 
the law of $B_{i+1}$ conditionally given $B_i=b$, is given by the
following transition. Let $a_1:=\supp(b)$, $b^*:=b-\delta_{a_1}$ and
$a_1^*:=\supp(b^*)$.
Let $(A,N)$ be distributed as $(A_1, \zeta'_{A_1})$. 
Then
\[
B_{i+1}:= \cases{ b^*+N\delta_A, &\quad if $A<
a_1^*$ and $A\neq a_1$,
\cr
b^*, &\quad otherwise. }
\]
%
\end{theorem}
\begin{pf}
Instead of the full information $(D_1(n), n\ge1)$, this sequence
starts with a single point measure
\[
B_1=D_1(A_1)\delta_{A_1},
\]
and at each step it proceeds by changing the weights of the existing
point masses and by adding at most one new point mass.

It is clear from the recursive definition of $(B_i;i\ge1)$ that for any
$i\ge1$ if \mbox{$b_{n,i}\neq0$}, then $b_{n,i}=D_i(n)$.
We first show that for any $i\ge1$ we have $a_{1,i}=A_i$ and
$b_{a_{1,i},i}=D_i(A_i)$. If we show that $b_{A_i,i}\neq0$, then,
since all other nonzero weights in $B_i$ satisfy $b_{n,i}=D_i(n)$, the
definition of $A_i$ will immediately imply that $a_{1,i}=A_i$ and
$b_{a_{1,i},i}=D_i(A_i)$. We do this by induction. The claim is clearly
true for $i=1$, so let us assume it is true for an arbitrary $i\ge1$.

Consider what the transition rule for $D_i$ tells us about the
relationship between $A_{i+1}$ and $A_i$. Recall for all $n>A_i$ we
have $D_{i+1}(n)=D_i(n)$, $D_{i+1}(A_i)=D_i(A_i)-1$, and for all
$n<A_i$ we have $D_{i+1}(n)=\zeta'_n$ from r.v.s drawn independently
of $D_i$. So,
\begin{eqnarray*}
A_{i+1}<A_i \quad&\Leftrightarrow&\quad \exists n<A_i
\mbox{ s.t. }\zeta_n'\neq0,
\\
A_{i+1}=A_i \quad&\Leftrightarrow&\quad D_i(A_i)>1
\mbox{ and } \zeta'_n=0,\forall n<A_i,
\\
A_{i+1}>A_i \quad&\Leftrightarrow&\quad D_i(A_i)=1
\mbox{ and } \zeta'_n=0, \forall n<A_i.
\end{eqnarray*}
In the first case, $A_{i+1}<A_i=a_{1,i}\le a_{1,i}^*$, so the point
mass at $A_{i+1}$ will be added to $B_{i+1}$ and $b_{A_{i+1},i+1}\neq
0$. In the second case, $b_{A_i,i+1}=D_{i+1}(A_i)=D_i(A_i)-1>0$, and
since $A_{i+1}=A_i$, we have $b_{A_{i+1},i+1}\neq0$. In the third
case, $b_{A_i,i+1}=D_{i+1}(A_i)=D_i(A_i)-1=0$, and also
$b_{a_{1,i},i}=D_i(A_i)=1$ implies $a_{1,i}^*=\min\{n>a_{1,i}\dvtx
b_{n,i}\neq0\}=\min\{n>A_i\dvtx  b_{n,i}\neq0\}$. Note that for all
$n>A_i$ for which $b_{n,i}\neq0$, we also have $b_{n,i+1}=b_{n,i}\neq
0$, and since $b_{A_i,i+1}=0$, the smallest value of $n$ with mass
$B_{i+1}(\{n\})>0$ before we potentially add a new mass is precisely at
$a_{1,i}^*=\min\{n>A_i\dvtx  b_{n,i}\neq0\}$. In this case
$a_{1,i}^*>a_{1,i}=A_i$ so we have
$D_{i+1}(a_{1,i}^*)=D_i(a_{1,i}^*)\neq0$, so by definition of
$A_{i+1}$ we must have $A_{i+1}\le a_{1,i}^*$. In case $A_{i+1}<
a_{1,i}^*$, the point mass at $A_{i+1}$ will be added to $B_{i+1}$ and
$b_{A_{i+1},i+1}\neq0$. In case $A_{i+1}= a_{1,i}^*$, no new mass is
added and the smallest of the nonzero masses in $B_{i+1}$ is at
$a_{1,i}^*=A_{i+1}$, as stated earlier, and again
$b_{A_{i+1},i+1}=b_{a_{1,i}^*,i+1}=b_{a_{1,i}^*,i}\neq0$. Hence, we
have shown by induction that $b_{A_{i},i}\neq0$ for all $i\ge1$, so
that $a_{1,i}=A_i$ and $b_{a_{1,i},i}= B(A_i,i)$.

Now consider the transition rule for $B_{i+1}$, conditionally given
$B_i$. For the already existing mass in $B_i$ the changes in weights
are given by the transition rule for $D_i$ to be
\begin{eqnarray*}
\sum_{n\ge1}D_{i+1}(n)\mathbf{1}_{\{b_{n,i}\neq0\}}
\delta_n &=& \bigl(D_i(A_i)-1 \bigr)
\delta_{A_i}+\sum_{n>A_i}D_i(n)
\mathbf{1}_{\{
b_{n,i}\neq0\}} \delta_n \\
&=& B_i -
\delta_{a_{1,i}}=B_i^*.
\end{eqnarray*}
We have an addition of a new point mass iff $A_{i+1}<a_{1,i}^*$. Since
$a_{1,i}^*\ge a_{1,i}=A_i$, this happens iff either:
\begin{longlist}[(ii)]
\item[(i)] $A_{i+1}<A_i$, or
\item[(ii)] $b_{a_{1,i},i}=1$ and
$A_i<A_{i+1}<a_{1,i}^*$.
\end{longlist}
The reason $A_i=A_{i+1}$ is not included in (ii) is that if
$b_{a_{1,i},i}=D_i(A_i)=1$ then $D_{i+1}(A_i)=D_i(A_i)-1=0$.
Let $\{\zeta'_n\}_{n\ge1}$ be a sequence of r.v.s drawn independently
from $D_i$. Then, by the transition rule for $D_i$:
\begin{longlist}[(ii)]
\item[(i)]holds iff $\exists n<A_i$ s.t. $\zeta_n'\neq0$, then,
$A_{i+1}:=\min\{n<A_i\dvtx \zeta_n'\neq0\}$;
\item[(ii)]holds iff $b_{a_{1,i},i}=1$, $\forall n<A_i$ $\zeta'_n=0$
and $\exists n<a_{1,i}^*$ s.t. $\zeta'_n\neq0$, then,
$A_{i+1}:=\min
\{A_i<n<a_{1,i}^*\dvtx \zeta_n'\neq0\}$.
\end{longlist}

In case\vspace*{1pt} (i) holds it is clear that the weight $D_{i+1}(A_{i+1})$ is
distributed as $\zeta'_{A_1}$ conditional on $A_1<A_i=a_{1,i}^*$ and
that this weight is independent of $(B_j; 1\le j\le i)$. We next argue
that this is also true in case (ii) holds. In this case
$A_i<A_{i+1}<a_{1,i}^*$ and since $D_{i+1}(A_{i+1})\neq0$ we must also
have $D_i(A_{i+1})\neq0$ because the transition rule for $D_{i+1}$
only allows zero entries to become nonzero for $n<A_i$. Moreover,
$D_i(A_{i+1})=D_{i+1}(A_{i+1})$ because all entries for $n>A_i$ remain
unchanged from step $i$ to $i+1$. Hence, we must have
$D_n(A_{i+1})=D_{i+1}(A_{i+1})$ for all $k<n\le i$ where
\[
k:=\max\{0\le j<i\dvtx A_j\ge a_{1,i}^*\}
\]
with the convention that $A_0:=+\infty$.\vadjust{\goodbreak}

We show that $A_k\neq A_{i+1}$ if $A_i<A_{i+1}<a_{1,i}^*$. If we had
$A_k=A_{i+1}$, then $b_{A_k,k}=D_k(A_k)\neq0$. The transition rule
for $D_{k+1}$ implies $D_{k+1}(A_k)=\break D_k(A_k)-1$. Since
$A_n<A_{i+1}=A_k$ for all $k+1\le n\le i$, iteratively applying the
transition rule for $D_{k+2},\ldots, D_{i+1}$ implies
$D_{k+2}(A_k)=D_{k+1}(A_k), \ldots,\break D_{i+1}(A_k)=D_i(A_k)$. Thus
\[
D_k(A_k)-1=D_{k+1}(A_k)=\cdots
=D_i(A_k)=D_{i+1}(A_k)=D_{i+1}(A_{i+1})
\neq0
\]
and $b_{A_k,k}=D_k(A_k)>1$. Then, by definition of the weights for
$B_{k+1},\ldots, B_{i+1}$, we have that $b_{A_k,k+1}=D_{k+1}(A_k)\neq0$
and the same weight iteratively remains as $b_{A_k,n}=D_n(A_k)$ for all
$k+1\le n\le i+1$. However, $b_{A_{i+1},i}=b_{A_k,i}\neq0$ contradicts
our assumption that $A_{i+1}<a_{1,i}^*$. 

Thus we must have $A_k>A_{i+1}$. Since $D_i(A_k)\neq0$ we also must
have $A_k\ge a_{1,i}^*$.
From the definition of weights for $B_i$ and $A_i<a_{1,i}^*$ we have
\[
b_{a_{1,i}^*,i}\neq0 \quad\Rightarrow\quad\exists k'< i \mbox{ s.t. }
a_{1,i}^*=A_{k'} \mbox{ and }A_n<A_{k'}
\mbox{ for all } k'<n\le i
\]
since the only point mass in a step $n$ not existing in the previous
step must be placed at $A_n$.
Now, $A_{k'}=a_{1,i}^*>A_{i+1}$ and the fact that we defined $k=\max\{
0\le j<i\dvtx A_j\ge A_{i+1}\}$ implies that $k=k'$, hence, $A_k= a_{1,i}^*$.

Since $A_k>A_{i+1}$, we must have $D_k(A_{i+1})=0$. By the same
argument as above of iteratively applying the transition rule for
$D_{k+2},\ldots, D_{i+1}$, the weights at $A_{i+1}$ satisfy
\[
D_{k+1}(A_{i+1})=\cdots=D_i(A_{i+1})=D_{i+1}(A_{i+1})
\neq0.
\]

Let us now consider the distribution of $D_{i+1}(A_{i+1})$ conditional
on the value of~$k$. Since $A_{i+1}<A_k$ and $D_{k+1}(A_{i+1})\neq0$
by the transition rule for $D_{k+1}$, the value of the weight
$D_{k+1}(A_{i+1})$ is a r.v. distributed as $\zeta'_{A_1}$ conditional
on $A_1<A_k$ and it is drawn independently of $(B_j; 1\le j\le k)$. So,
$D_{i+1}(A_{i+1})=D_{k+1}(A_{i+1})$ is distributed as $\zeta'_{A_1}$
conditional on $A_1<A_k=a_{1,i}^*$ and is independent of $(B_j; 1\le
j\le k)$.
Since $A_n<A_{i+1}$ for all $k< n\le i$, the transition rule for
$D_{k+1},\ldots, D_i$ implies that all of the subsequently added point
masses in $B_{k+1},\ldots, B_i$ are independent of $D_{k+1}(A_{i+1})$,
so the value of $D_{i+1}(A_{i+1})$ is independent of $(B_j; 1\le j\le
i)$. Finally, integrating over $k$, we have that $D_{i+1}(A_{i+1})$ is
distributed as $\zeta'_{A_1}$ conditional on $A_1<a_{1,i}^*$ and is
independent of $(B_j; 1\le j\le i)$. 

We have now shown that when a new point mass is added at $A_{i+1}$,
either $A_{i+1}:=\min\{n<A_i\dvtx \zeta_n'\neq0\}$, or $b_{a_{1,i},i}=1$
and $A_{i+1}:=\min\{A_i<n<a_{1,i}^*\dvtx \zeta_n'\neq0\}$, where the
sequence $\{\zeta'_n\}_{n\ge1}$ are r.v.s drawn independently from
$(B_j; j\le i)$. In either case, the weight of the new point mass at
$A_{i+1}$ is distributed as $\zeta'_{A_1}$ conditional on $A_1<a_{1,i}^*$.

Let $A:=\min\{n\ge1\dvtx \zeta'_n\neq0\}$ and $N:=\zeta'_A$. Since $\{
\zeta'_n\}_{n\ge1}$ are independent of $(B_j; 1\le j\le i)$, so are $A$
and $N$, and by Theorem~\ref{thmcpp}, $(A,N)$ is distributed as
$(A_1,\zeta'_{A_1})$. Conditional on the value of $B_i$, we have that
(i) holds iff $A<A_i$, while (ii) holds iff $b_{a_{1,i},i}=1$ and
$A_i<A<a_{1,i}^*$. Putting (i) and (ii) together\vadjust{\goodbreak} with the definition of
$a_{1,i}^*$ we have, conditionally on $(B_j;1\le j \le i)$, an addition
of a new point mass iff $A<a_{1,i}^*$ and $A\neq a_{1,i}$. Moreover,
the weight of the newly added point mass is distributed as $\zeta
'_{A_1}$ conditional on $A_1<a_{1,i}^*$, or equivalently it is
distributed as $N$ conditional on $A<a_{1,i}^*$. We also showed that
given $(B_j;1\le j \le i)$, the point masses existing in $B_i$ change
in the next step to produce a re-weighted point mass measure equal to
$B_i-\delta_{a_1}$. Altogether, given $(B_j;1\le j \le i)$, the next
step of the sequence depends on $B_i$ only, with the transition rule
that if $A<a_{1,i}^*$ and $A\neq a_{1,i}$ in the next step we have
$B_{i+1}=B_i-\delta_{a_{1,i}}+N\delta_A$, and otherwise in the next
step we have $B_{i+1}=B_i-\delta_{a_{1,i}}$.\vspace*{-2pt}
\end{pf}

In the course of the above proof, we have also shown that if $t_1:=\min
\{t\ge1\dvtx A_t>A_1\}$ and $N_1:=\#\{1\le j<t_1\dvtx  A_j=A_1\}$, then
$N_1=D_1(A_1)$, because in order to add the first new point mass at a
level greater than $A_1$, we must in the course of $1\le j<t_1$ have
exactly enough steps at which $A_j=A_1$ that will exhaust all of the
weight $D_1(A_1)$ of the point mass at $A_1$.

Analogously for each $i\ge1$, if we let
\[
t_i:=\min\{t\ge i\dvtx A_t>A_i\}\quad \mbox{and}\quad
N_i:=\#\{i\le j<t_i\dvtx  A_j=A_i
\},
\]
then $D_i(A_i)=N_i$.
If, furthermore, for each $i\ge1$ and $k\ge i$ we define
\[
N_{ik}:=\#\{k\le j<t_i\dvtx A_j=A_i
\},
\]
then by a similar argument for $k\ge i$ we have $N_{ik}=D_k(A_i)$. Note
that $N_{ii}\equiv N_i$.
Then, for the sequence of finite point measures we have that for all
$k\ge1$,
%
\begin{equation}
\label{eqaltcoalesc} B_k=\sum_{i\le k}N_{ik}
\delta_{A_i}\mathbf{1}_{\{A_i< A_{i'},
\forall
{i'}<i\dvtx  N_{{i'}k}\neq0\}}.
\end{equation}
It is easily checked that $B_{k+1}$ correctly updates the weight of
each existing point mass from $B_k$, and allows a new point mass to be
added only if it is in a location smaller than all mass existing in
$B_k$ whose weights in $B_{k+1}$ remain nonzero.
We will see this formula again when we discuss the coalescent process
in the continuous case.\vspace*{-2pt}
%
\begin{rem} In the discrete case the coalescence times $A_i$ take on
integer values which may occur again after they have appeared for the
last time $t_i$ in a subtree. In the continuous case, analogously
defined coalescence times have a law that is absolutely continuous
w.r.t. Lebesgue measure, their values a.s. never occur again after they
have appeared for the last time in a subtree, hence, in the continuous
case there is no need to use separation times $t_i$ in the definition
of counters $N_i$ and~$N_{ik}$.\vspace*{-2pt}
\end{rem}

\section{From discrete to continuous: The great-aunt measure}\vspace*{-2pt}

\subsection{Definition of the great-aunt measure}

Now that we have a simple description of the coalescent point process
with multiplicities for a BGW branching\vadjust{\goodbreak} process, our goal is to do the
same for the continuous-state branching process. In the discrete case
we used the process $(D_i;i\ge1)$ to describe the number of surviving
younger offshoots of ancestors of an individual as a sequence indexed
by generations backward in time. Since in the CSB case the standing
population is not discrete, we cannot use a process indexed by the
standing population. Instead, we use a spine decomposition of the
lineage of a surviving individual, which will record the level (i.e.,
generation in the discrete, and height in the continuous case) and the
number of all offshoot subtrees in the individual's genealogy. We first
provide the law of the spine decomposition of the first individual in
the standing population in the BGW case and relate it to our previous
results. At this point we would like to emphasize that results for the
spine decomposition of BGW process are not new (see references at the
end of the subsection), and that we make use of the spine decomposition
here only as a tool that will enable us to describe the analogue of the
coalescent point process with multiplicities for the CSB process later.

We give some new definitions to describe the spine decomposition of
$(0,1)$, the first individual in generation $0$, of a BGW process
within its monotone planar embedding. For $n\ge1$, we denote by $\phi
_n$ the set of great-aunts of $(0,1)$ at generation $-n+1$, that is,
\[
\phi_n:=\bigl\{\mbox{daughters of } \bigl(-n,\mathfrak{a}_1(n)
\bigr) \mbox{ excluding }\bigl(-n+1,\mathfrak{a}_1(n-1)\bigr)\bigr\},
\]
and by $\varphi_n$ its cardinal, $\varphi_n:=\# \phi_n$. In other
words, $\varphi_n$ is the number of offshoots of $(-n,\mathfrak
{a}_1(n))$ not counting the lineage of $(0, 1)$ itself.
The set $\phi_n$ can be divided into
\[
\alpha_n:=\bigl\{\mbox{sisters of } \bigl(-n+1,\mathfrak{a}_1(n-1)
\bigr) \mbox{ with labels } (-n+1,k) \mbox{ s.t. } k<\mathfrak{a}_1(n-1)
\bigr\}
\]
the older offshoots of $(-n,\mathfrak{a}_1(n))$, and
\[
\rho_n:= \bigl\{\mbox{sisters of } \bigl(-n+1,\mathfrak{a}_1(n-1)
\bigr) \mbox{ with labels } (-n+1,k) \mbox{ s.t. } k>\mathfrak{a}_1(n-1)
\bigr\}
\]
the younger offshoots of $(-n,\mathfrak{a}_1(n))$.
We call the sequence $(\rho_n;n\ge1)$ the \textit{great-aunt measure}.
We let $(Z_k^{\alpha_n}; 0\le k\le n-1)$ and $(Z_k^{\rho_n}; 0\le
k\le
n-1)$ be the processes counting the descendants of those $\alpha_n$ and
$\rho_n$ individuals, respectively, at successive generations $-n+1,\ldots, 0$.

An important observation now is that we do not need the whole infinite
embedding of trees to make the previous definitions. The descending
tree of $(-h,\mathfrak{a}_1(h))$ is a planar BGW tree conditioned to
have alive individuals at generation $h$, with the lineage of
individual $(0,1)$ being the left-most lineage with alive descendants
at generation $h$. Consequently, provided that we only consider indices
$n\le h$, our definitions make sense for any conditioned BGW planar
tree. A number of results have already been proved in~\cite{Gei2} for
spine decomposition of BGW process conditioned on survival at a given
generation, and we make use of them here. The infinite embedding of
trees that we introduced extends these\vadjust{\goodbreak} results in a way, as we are
considering an arbitrarily large standing population that may not [in
(sub)critical case] descend from the same ancestor $n$ generations back
in the past. Considering trees whose roots are individuals
$(-n,\mathfrak{a}_1(n)), (-n, \mathfrak{a}_1(n)+1), \ldots\,$, it is easy
to see that they are independent and identically distributed as ${\cal
T}(n,1)$, the tree descending from $(-n,\mathfrak{a}_1(n))$.

We recall some standard notation commonly used with linearly ordered
planar trees.
The \textit{Ulam--Harris--Neveu} labeling of a planar tree assumes that
each individual (vertex) of the tree is labeled by a finite word $u$ of
positive integers whose length $\vert u\vert$ is the generation, or
height, of the individual. A~rooted, linearly ordered planar tree
${\cal T}$ is a subset of the set ${\cal U }$ of finite words of integers
\[
{\cal U }= \bigcup_{n\ge0} \NN^n,
\]
where $\NN^0$ is the empty set.
More specifically, the root of $\cal T$ is denoted $\varnothing$ and
her offspring are labeled $1,2, 3,\ldots$ from left to right. The
recursive rule is that the offspring of any individual labeled $u$ are
labeled $u1, u2, u3,\ldots$ from left to right, where $ui$ is the mere
concatenation of the word $u$ and the number $i$.
The \textit{depth-first search} is the lexicographical order associated
with Ulam--Harris--Neveu labeling (see Figure~\ref{figrandomwalk},
where the depth-first search gives the order $\varnothing, 1, 11, 12,
2, 211, 212, 3, 31, 311, 32, \ldots$).

Now fix $h\in\NN$ and assume that generation $h$ is nonempty. Let $x_h$
be the Ulam--Harris--Neveu label of the first individual in depth-first
search with height $h$. We denote by $(x_h\mid-n)$ the ancestor of
$x_h$ at generation $h-n$. Then for $1\le n\le h$, we can define
$\alpha'_n$ and $\rho'_n$ to be the number of daughters of $(x_h\mid-n)$
ranked smaller and larger, respectively, than $(x_h\mid-n+1)$, and
define $\varphi'_n:=\alpha'_n+\rho'_n$. We can also let $Z^{\alpha'_n}$
and $Z^{\rho'_n}$ be the processes counting the descendants of those
$\alpha'_n$ and $\rho'_n$ individuals, respectively.
From now on, we assume that the tree $\cal T$ has the law of a planar
BGW tree with offspring distributed as $\xi$ and conditioned to have
alive individuals at generation $h$. Then it is easily seen that
$(\alpha'_n, \rho'_n, Z^{\alpha'_n}\circ\kil_n, Z^{\rho'_n}\circ
\kil_n;1\le n\le h)$ (where $\kil_n$ means killing at time $n$) has the
same distribution as $(\alpha_n, \rho_n, Z^{\alpha_n}, Z^{\rho
_n};1\le
n\le h)$, so from now on we remove primes. The following result
provides the joint law of these random variables, which was already
shown in~\cite{Gei2}, and the results below are just a restatement of
Lemma 2.1 from~\cite{Gei2}. Recall that $p_n=\PP(Z_n\not=0\mid Z_0=1)$.
%
\begin{prop}
\label{propgreataunts}
Conditional on the values of $(\alpha_n, \rho_n; n\ge1)$, the
processes $(Z^{\alpha_n}, Z^{\rho_n};n\ge1)$ are all independent,
$Z^{\rho_n}$ is a copy of $Z$ started at $\rho_n$ and killed at time
$n$ and $Z^{\alpha_n}$ is a copy of $Z$ started at $\alpha_n$,
conditioned on being zero at time $n-1$. In addition, the pairs
$(\alpha_n, \rho_n; n\ge1)$ are independent and distributed as follows:
%
\[
\PP(\alpha_n = j, \varphi_n=k) = \PP(\xi=k+1)
\frac
{p_{n-1}}{p_{n}}(1-p_{n-1})^{j},\qquad k\ge j\ge0.
\]
%
\end{prop}
%
%
\begin{rem}
\label{rembinom}
Recall that $D_1(n)$ is the number of sisters of $(-n+1, \mathfrak
{a}_1(n-1))$ with alive descendants in the standing population.
Thus, an immediate corollary of Proposition~\ref{propgreataunts}
is that the random variables $(D_1(n);n\ge1)$ are independent and
that, conditionally on $\rho_n$, $D_1(n)$ is a binomial r.v. with
parameters $\rho_n$ and $p_{n-1}$. Since, according to Theorem
\ref{thmcpp}, $D_1(n)$ is distributed as $\zeta'_n$, we also have
\[
\EE\bigl(s^{\zeta'_n} \bigr)=\EE\bigl( (1-p_{n-1}+p_{n-1}s)^{\rho
_n}
\bigr)
\]
as one can easily check.
\end{rem}

Note that the last statement in Proposition~\ref{propgreataunts}
can be viewed as an inhomogeneous spine decomposition of the
``\textit{ascendance}'' of the surviving particles in conditioned BGW trees.
More standard spine decompositions are well known for the
``\textit{descendance}'' of conditioned branching trees (see, e.g., \cite
{Ev,L-PTRF,L-EJP,LPP,Gei1,CRW,K}) (the idea of spine decompositions
originating from~\cite{K}). To bridge the gap between both aspects,
notice that in the (sub)critical case, $\PP(\varphi_n+1 = k)= \PP
(\xi
=k) (1-(1-p_{n-1})^k)/p_n$ converges to $k\PP(\xi=k)/\EE(\xi)$ as
$n\to
\infty$, which is the size-biased distribution of $\xi$. This
distribution is known to be the law of offspring sizes on the spine
when conditioning the tree on infinite survival.

\subsection{A random walk representation}
\label{subsecrwr}

We next show how to recover any truncation $(\rho_n;n\le h)$ of the
great-aunt measure from the paths of a conditioned, killed random walk.
This will be particularly useful when we define the analogue of the
great-aunt measure for continuous-state branching processes in the next
subsection.
Our result makes use of a well-known correspondence between a BGW tree
and a downward-skip-free random walk introduced and studied in~\cite{BK,LGLJ}.

Let us go back to the planar tree $\cal T$.
We denote by $v_n$ the word labeling the $n$th individual of the tree
in the depth-first search. For any integers $i>j$ and any finite word
$u$, we say that $ui$ is ``younger'' than $uj$. Also, we will write
$u\prec v$ if $u$ is an ancestor of $v$, that is, there is a sequence
$w$ such that $v=uw$ (in particular, $u\prec u$). Last, for any
individual $u$ of the tree, we let $r(u)$ denote the number of younger
sisters of $u$, and for any integer $n$, we let $W_n:=0$ if
$v_n=\varnothing$, and if $v_n$ is any individual of the tree different
from the root, we let
\[
W_n:=\sum_{u\prec v_n} r(u).
\]
The height $H_n$, or generation, of the individual visited at time $n$
can be recovered from $W$ as
%
\begin{equation}
\label{eqndefdiscreteheight} H_n:= \vert v_n \vert=
\#\Bigl\{1\le k\le n\dvtx  W_k = \min_{k\le j \le n}W_j\Bigr
\}.
\end{equation}
See Figure~\ref{figrandomwalk} where $(W_n;n=1,\ldots,15)$ is given
until visit of $v_{15}=3221$ whose height is $H_{15}=4$.

In the case when $\cal T$ is a BGW tree with offspring distributed as
$\xi$, then it is known~\cite{BK,LGLJ} that the process $(W_n;n\ge1)$
is a random walk started at 0, killed upon hitting $-1$, with steps in
$\{-1,0,1,2,\ldots\}$ distributed as $\xi-1$.
%

Fix $h\in\NN$ and set
\[
\sigma_h:=\min\{n\ge1\dvtx  H_n = h\}.
\]
Writing $T_j$ for the first hitting time of $j$, in particular,
$T_{-1}$ for the first hitting time of $-1$, we get
%
\begin{equation}
\label{eqnmaxdiscreteheight1} \Bigl\{\max_{0\le j \le T_{-1}} H_j
< h\Bigr\} = \{Z_h=0\} = \{\sigma_h \ge T_{-1}
\}.
\end{equation}

Now, when $\sigma_h<T_{-1}$, we let $I_{\sigma_h}$ denote the future
infimum process of the random walk killed at $\sigma_h$
\[
I_{\sigma_h}^j:= \min_{j\le r\le\sigma_h} W_r,
\]
and we let $t_h=0,\ldots,t_0=\sigma_h$ denote the successive record
times of $I_{\sigma_h}$
\[
t_k:=\max\bigl\{j<t_{k-1}\dvtx  W_j =
I_{\sigma_h}^j\bigr\},\qquad 1\le k \le h.
\]
Also observe that by definition of $\sigma_h$, we must have $t_1=
t_0-1$. Last, we use the notation $\Delta_j W:=W_{j+1}-W_j$. A
straightforward consequence of~\cite{BK,LGLJ} is that when generation
$h$ is nonempty, $t_k$ is the visit time of $(x_h\mid-k)$, that is,
the unique integer $n$ such that $v_n=(x_h\mid-k)$ (where $x_h$ is the
first individual in depth-first search with height $h$). Furthermore,
\[
\varphi_k=\Delta_{t_k} W=W_{t_k+1}-W_{t_k}
\]
and
%
\[
\alpha_k=W_{t_k+1}- \min_{t_k+1\le j \le t_0}W_j
\quad\mbox{and}\quad \rho_k= \min_{t_k+1\le j \le t_0}W_j -
W_{t_k}.
\]
In particular, we can check ${\alpha}_1=0$, and ${\varphi}_k={\alpha
}_k+{\rho}_k$. Now note that, by definition of $t_h,\ldots,t_0$, we have
$I^{t_k}_{\sigma_h}= W_{t_k}$ and for all $j\dvtx  t_k<j<t_{k-1}$ we have
$I^{j}_{\sigma_h}=W_{t_{k-1}}$, so
\[
{\alpha}_k=W_{t_k+1}-I^{t_k+1}_{\sigma_h}=W_{t_k+1}-W_{t_{k-1}}
\quad\mbox{and}\quad {\rho}_k=I^{t_k+1}_{\sigma_h}-W_{t_k}=W_{t_{k-1}}-W_{t_k}.
\]
An illustration of these claims can be seen in Figure~\ref{figrandomwalk}.


Last, observe that
\begin{eqnarray*}
\sum_{k=1}^h {\rho}_k f(k)&=&
\sum_{k=1}^h \bigl(I_{\sigma_h}^{t_k+1}
- W_{t_k}\bigr) f(k)=\sum_{k=1}^h
\sum_{j=0}^{\sigma_h-1}\indicbis{j=t_k}
\bigl(I_{\sigma_h}^{j+1} - W_{j}\bigr)
f(h-H_j)\\
&=&\sum_{j=0}^{\sigma_h-1}
\Delta_j I_{\sigma_h}^jf(h-H_j)
\end{eqnarray*}
since, if $j=t_k$ for some $k$, then
\[
\Delta_j I^j_{\sigma_h}=\Delta_{t_k}
I^{t_k}_{\sigma
_h}=I^{t_k+1}_{\sigma_h}-I^{t_k}_{\sigma_h}=I^{t_k+1}_{\sigma_h}-W_{t_k},
\]
otherwise, $t_k<j<t_{k-1}$ for some $k$, and
\[
\Delta_{j} I^{j}_{\sigma_h}=I^{j+1}_{\sigma_h}-I^{j}_{\sigma
_h}=W_{t_{k-1}}-W_{t_{k-1}}=0.
\]
This is recorded in the following statement, which we will use later as
a distributional equality for a random walk $W$ conditioned on $\max
_{1\le j \le T_{-1}} H_j\ge h$.

\begin{prop}
\label{proprwgreat-aunts}
For any $f\dvtx \NN\longrightarrow\RR_+$ with bounded support,\break $h>
\operatorname{Supp}(f)$, if $\max_{1\le j \le T_{-1}} H_j\ge h$, then
%
\[
\langle\rho, f\rangle:= \sum_{k\ge1}
\rho_k f(k)=\sum_{j=0}^{\sigma_h-1}
\Delta_j I_{\sigma_h}^jf(h-H_j),
\]
%
where $\Delta_j I_{\sigma_h}^j=I_{\sigma_h}^{j+1}-I_{\sigma_h}^j$.
\end{prop}
\subsection{A continuous version of the great-aunt measure}

In Section~\ref{sec2}, we gave a consistent way of embedding trees with an
arbitrary size of the standing population, each descending from an
arbitrarily old founding ancestor, so that the descending subtree of
each vertex is a BGW tree. The natural analog of this presentation is
the flow of subordinators introduced by Bertoin and Le Gall in \cite
{BLG}. Because the Poissonian construction of this flow displayed in
\cite{BLG3}, Section 2, only holds for critical CSB processes and
without Gaussian component, and because it is rather awkward to use in
order to handle the questions we address here, we will now define an
analogue of the great-aunt measure, using the genealogy of
continuous-state branching processes introduced in~\cite{LGLJ} and
further investigated in~\cite{DLG}.

We start with a L\'{e}vy process $X$ with no negative jumps and Laplace
exponent $\psi$ started at $x>0$, such that $\psi'(0+)\ge0$, so that
$X$ hits 0 a.s. As specified in \mbox{\cite{L1,L2}}, the path of $X$ killed
upon reaching 0 can be seen as the (jumping) contour process of a
continuous tree whose ancestor is the interval $(0,x)$. For example,
the excursions of $X$ above its past infimum $I$ are the contour
processes of the offspring subtrees of the ancestor. Almost surely for
all $t$, it is possible to define the \textit{height} (i.e., generation)
$H_t$ of the point of the tree visited at time $t$, as
%
\begin{equation}
\label{eqndefheightprocess} H_t:= \lim_{k\to\infty}
\frac{1}{\vareps_k}\int_0^t \indic
{X_s<I_t^s+\vareps_k} \,ds <\infty,
\end{equation}
where $(\vareps_k)$ is some specified vanishing positive sequence, and
\[
I_t^s:= \inf_{s\le r\le t} X_r.
\]
With this definition, one can recover the population size at generation
$a$ as the density $Z_a$ of the occupation measure of $H$ at $a$, that
is, the (total) local time of $H$ at\vadjust{\goodbreak} level $a$. It is proved in \cite
{LGLJ,DLG} that this local time exists a.s. for all $a$ and that
$(Z_a;a\ge0)$ is a continuous-state branching process with branching
mechanism $\psi$.

From now on, we will deal with a general branching mechanism $\psi$
characterized by its L\'{e}vy--Khinchin representation
\[
\psi(\lbd) = a\lbd+\beta\lbd^2 +\int_{(0,\infty)}
\Lambda(dr) \bigl(e^{-\lbd
r} -1 +\lbd r \indic{r<1}\bigr),
\]
where $\beta>0$ is the \textit{Gaussian component} and $\Lambda$ is a
positive measure on $(0,\infty)$ such that $\int_{(0,\infty
)}(1\wedge
r^2)\Lambda(dr)<\infty$, called the \textit{L\'{e}vy measure}.
We will denote by $X$ a L\'{e}vy process with Laplace exponent $\psi$
(started at 0 unless otherwise stated), and by $Z$ a continuous-state
branching process, or CSB process, with branching mechanism $\psi$.

We (only) make the following two assumptions. First,
\[
\sup\bigl\{\lbd\dvtx  \psi(\lbd)\le0\bigr\}=:\eta<\infty,
\]
so that $X$ is not a subordinator (i.e., it is not a.s. nondecreasing). Second,
\[
\int_{[\eta+1,\infty)} \frac{du}{\psi(u)}<\infty,
\]
so that $Z$ either is absorbed at 0 in finite time or goes to $\infty$,
and $H$ has \textit{a.s. continuous sample paths}~\cite{LGLJ}. This also
forces $\int_{(0,1)}r\Lambda(dr)$ to be infinite, so that the paths of
$X$ have \textit{infinite variation}.
We further set
\[
\phi(\lbd):=\int_{[\lbd,\infty)}\frac{du}{\psi(u)},\qquad \lbd>\eta,
\]
and $v$ the inverse of $\phi$ on $(\eta,+\infty)$
\[
v(x):=\phi^{-1}(x),\qquad x\in(0,\infty).
\]
Notice that $v$ is nonincreasing and has limit $\eta$ at $+\infty$. It
is well known (e.g.,~\cite{L1}) that if $Z$ is started at $x$, then it
is absorbed at 0 before time $t$ with probability $e^{-xv(t)}$. Also,
if $N$ denotes the excursion measure of $X-I$ away from 0 under $\PP$
(normalized so that $-I$ is the local time), then by~\cite{DLG},
Corollary 1.4.2,
\[
N(\sup H >a) = v(a),\qquad a>0.
\]

Here, we want to allow the heights of the tree to take negative values.
To do this, we start with a measure which embodies the mass
distribution broken down on heights, of the population whose
descendances have not yet been visited, in the same vein as in \cite
{LGLJ,DLG}, but with negative heights. In the usual setting, the mass
distribution $\rho_t$ of the population whose descendances have not yet
been visited by the contour process $X$ by time $t$, is defined by
\[
\langle\rho_t,f\rangle=\int_{[0,t]}d_sI_{t}^sf(H_s)
\]
for any nonnegative function $f$, where on the right-hand side we mean
integrating the function $s\mapsto f(H_s)$ with respect to the
Stieltjes measure associated with the function $s\mapsto I_{t}^s$. Then
$\rho_t([a,b])$ is the mass of the tree between heights $a$ and $b$
whose descendants have not yet been visited by the contour process.

Here, we will start with a random positive measure $\rho^0$ on
$[0,+\infty)$, with the interpretation that $\rho^0([a,b])$ is the mass
of the tree between heights $-b$ and $-a$ whose descendants have not
yet been visited by the contour process. This measure is the exact
analogue of the great-aunt measure of the previous subsection.
%
\begin{dfn}
\label{dfnrho}
For every $t>0$, set $\pi^{(t)} (dr):= p(t,r) \,dr$, where
\[
p(t,r):=e^{v(t)r} \int_{(r,\infty)} e^{-v(t)z}
\Lambda(dz),\qquad r>0.
\]
We define $\rho^0$ in law by
\[
\rho^0(dx):= \beta \,dx+ \sum_{t\dvtx \Delta_t>0}
\Delta_t\delta_t(dx),
\]
where $\delta$ denotes a Dirac measure, and $(t,\Delta_t)$ is a Poisson
point measure with intensity measure $dt\, \pi^{(t)}(dr)$.
\end{dfn}

\subsection{Convergence of the great-aunt measure}
\label{subsecconvofg-a}

We can now prove a theorem yielding two justifications for the
definition of the measure $\rho^0$. First, we show that $\rho^0(dx)$ as
defined above is indeed the mass, at height $h-x$, of the part of the
tree whose descendants have not yet been visited, either by a
long-lived contour process ($h\to\infty$) or under the measure
$N(\cdot
\mid\sup H > h)$. Second, we show the convergence of the appropriately
re-scaled discrete great-aunt measures to the measure $\rho^0$, as the
BGW processes approach the CSB process $Z$ with branching mechanism
$\psi$. In the next section we will show how $\rho^0$ allows us to
define the coalescent point process with multiplicities for the CSB
process, and help us establish convergence from the appropriately
re-scaled point process $(B_i;i\ge1)$.

We assume there exists a sequence $(\gamma_p; p\ge1)$, $\gamma_p\to
\infty$ as $p\to\infty$, and a sequence of random variables $(\xi_p;p\ge
1)$, such that, if $W^{(p)}$ denotes the random walk with steps
distributed as $\xi_p -1$, the random variables
$(p^{-1}W^{(p)}_{[p\gamma_p]})$ converge in law to $X_1$, where $X$
denotes the L\'{e}vy process with Laplace exponent $\psi$. Then it is
known (\cite{Gri}, Theorems 3.1 and 3.4) that if $Z^{(p)}$ denotes the
BGW process started at $[px]$ with offspring size distributed as $\xi
_p$, then the re-scaled processes $(p^{-1}Z^{(p)}_{[\gamma_p t]};t\ge
0)$ converge weakly in law in Skorokhod space to the CSB process $Z$
with branching mechanism $\psi$, started at $x$.

In the following statement, we denote by $\rho^{(p)}$ the great-aunt
measure associated to the offspring size $\xi_p$. We have to make the
following\vadjust{\goodbreak} additional assumptions: if $f^{(p)}$ denotes the p.g.f. of
$\xi_p$ and $f^{(p)}_n$ its $n$th iterate, then for each $\delta>0$,
\[
\mathop{\lim\inf}_{p\to\infty}f^{(p)}_{[\delta\gamma_p]}(0)^p>0.
\]

Convergence results in (iii) below rely heavily on results already
established by Duquesne and Le Gall in~\cite{DLG} on convergence of
appropriately re-scaled random walks $W^{(p)}$ and their height
processes $H^{(p)}$ to the L\'{e}vy process $X$ and its height process
$H$. In particular, technical conditions such as the one above are
justified in~\cite{DLG}, Section 2.3.
%
\begin{theorem}
\label{thmpropertiesrho}
Let $\sigma_h$ denote the first time that $H$ hits $h>0$.
\begin{longlist}[(iii)]
\item[(i)] For any $h>0$ and for any nonnegative Borel function
$f$ vanishing outside $[0,h]$, the random variable $\int_{[0,{\sigma
_h}]}d_sI_{\sigma_h}^sf(h-H_s)$, under $N(\cdot\mid\sup H >h)$, has
the same law as $\langle\rho^0,f\rangle$;
\item[(ii)] For any nonnegative Borel function $f$ with compact
support, as $h\to\infty$ the random variables $\int_{[0,{\sigma
_h}]}d_sI_{\sigma_h}^sf(h-H_s)$ under $\PP(\cdot\mid\sigma_h
<\infty
)$ converge in distribution to $\langle\rho^0,f\rangle$;
\item[(iii)] For any nonnegative Borel function $f$ with compact
support and sequence of nonnegative continuous functions $\{f_p\}$
such that $f_{p,\gamma_p}(x_p):=f_p(\gamma_p^{-1} x_p)\to f(x)$
whenever $\gamma_p^{-1}x_p\to x$, 
the random variables $p^{-1} \langle\rho^{(p)},f_{p,\gamma_p}\rangle$
converge in
distribution to $\langle\rho^0,f\rangle$ as $p\to\infty$.
\end{longlist}
\end{theorem}
\begin{pf}
Let us prove (i). It is known~\cite{DLG} that a.s. for all $t$ the
inverse of the local time on $[0,t]$ of the set of increase times of
$(I_t^s;s\in[0,t])$ has drift coefficient $\beta$, so that
\[
\int_{[0,t]}d_sI_{t}^s
f(H_s) = \beta\int_0^{H_t} f(x) \,dx +
\sum_{s\in
[0,t]} f(H_s) \Delta_s
I_t^s,
\]
where the sum is taken over all times $s$ at which $(I_t^s;s\in[0,t])$
has a jump, whose size is then denoted $\Delta_s I_t^s$. This set of
times will be denoted ${\cal J}_h$ when $t=\sigma_h$. As a consequence,
it only remains to prove that the random point measure $M_h$ on
$(0,h)\times(0,\infty)$ defined by
\[
M_h:=\sum_{s\in{\cal J}_h} \delta
\bigl(H_s, \Delta_sI_{\sigma_h}^s\bigr),
\]
where the sum is zero when ${\cal J}_h$ is empty ($\sigma_h=\infty$),
is an inhomogeneous Poisson point measure with the correct intensity
measure. More precisely, we are going to check that for any
nonnegative two-variable Borel function~$f$
\begin{eqnarray*}
N \biggl(\int_{(0,h)\times(0,\infty)} M_h(dt \,dr) f(h-t,r)\Bigm| \sup H
>h \biggr)
= \int_0^h dt \int
_0^\infty dr \,p(t,r)f(t,r),
\end{eqnarray*}
that is,
\[
N\int_{(0,h)\times(0,\infty)} M_h(dt \,dr) f(h-t,r) = v(h) \int
_0^h dt \int_0^\infty
dr\, p(t,r)f(t,r).
\]

First notice that
\begin{eqnarray*}
&&
N\int_{(0,h)\times(0,\infty)} M_h(dt \,dr) f(h-t,r) \\
&&\qquad= N\sum
_{s\dvtx\Delta X_s
>0} \indicbis{s<\sigma_h}f\bigl(h-H_s,
\Delta X_{s}+K_s'\bigr)\indicbis{\Delta
X_{s} >-K_s' },
\end{eqnarray*}
where $K_s'$ is the global infimum of the shifted path $X'$
\[
X'_u:= X_{s+u}-X_s,\qquad 0\le u \le
\sigma_h-s.
\]
But\vspace*{1pt} on $\{s<\sigma_h\}$, $\sigma_h-s$ is also $\sigma'_{h-H_s}$ (with
obvious notation), so by predictable projection and by the compensation
formula applied to the Poisson point process of jumps, we get
\[
N\int_{(0,h)\times(0,\infty)} M_h(dt \,dr) f(h-t,r) = N\int
_0^{\sigma_h} ds \int_{(0,\infty)}
\Lambda(dz) G_f(z,H_s),
\]
where (with the notation $I$ for the current infimum)
\[
G_f(z,t) =\EE_{0} \bigl[f(h-t, z+I_{\sigma_{h-t}})
\indicbis{-I_{\sigma
_{h-t}}<z} \bigr].
\]
Now $-I_{\sigma_{h-t}}$ is the local time at the first excursion of
$X-I$ with height larger than $h-t$ so it is exponentially distributed
with parameter $N(\sup H> h-t)= v(h-t)$. As a consequence,
\[
G_f(z,t) =v(h-t)\int_0^z dr
f(h-t,r)e^{-(z-r)v(h-t)},
\]
which yields
\[
N\int_{(0,h)\times(0,\infty)} M_h(dt \,dr) f(h-t,r) = N\int
_0^{\sigma
_h}ds\, F_f(h-H_s)
\]
with
\begin{eqnarray*}
F_f(t) &= &\int_{(0,\infty)} \Lambda(dz)
G_f(z,h-t)
\\
&=& v(t) \int_{(0,\infty)} \Lambda(dz) \int_0^z
dr \,f(t,r) e^{-(z-r)v(t)}
\\
&=& v(t) \int_0^\infty dr\, f(t,r)
e^{rv(t)}\int_{(r,\infty)}\Lambda(dz) e^{-zv(t)}
\\
&=& v(t)\int_0^\infty dr\, p(t,r) f(t,r).
\end{eqnarray*}
So we only have to verify that for any nonnegative Borel function $g$,
\[
N\int_0^{\sigma_h}ds\,g(h-H_s)= \int
_0^h dt\,\frac{v(h)}{v(t)} g(t).
\]
Due to results in~\cite{LGLJ}, there indeed is a jointly measurable
process $(Z(a,t); \break a,t\ge0)$ such that a.s. for all $a$,
\[
\int_0^{a}ds\,g(H_s)= \int
_0^{\infty} dt\,Z(a,t) g(t).
\]
In particular,
\[
N\int_0^{\sigma_h}ds\,g(h-H_s)= \int
_0^{h} dt\,N\bigl(Z(\sigma_h,t)\bigr)
g(h-t),
\]
so we just need to check that $N(Z(\sigma_h,t))={v(h)}/{v(h-t)}$, that is,
\[
N\bigl(Z(\sigma_h,t)\mid\sup H >h\bigr)=\frac{1}{v(h-t)}.
\]
But, conditional on $\sup H >h$, $\sigma_t <\infty$ and $Z(\sigma_h,t)$
is the local time of $H$ at level~$t$ between $\sigma_t$ and $\sigma
_h$, which is exponential with parameter $N(\sup H>h-t)= v(h-t)$,
hence, the claimed expectation.

We proceed with (ii), which readily follows from (i). Indeed, for any
$h$ larger than $\sup\{x\dvtx  f(x)\not=0\}$,
\[
\int_{[0,{\sigma_h}]}d_sI_{\sigma_h}^sf(h-H_s)=
\int_{[0,{\sigma
_h}]}d_sI_{\sigma_h}^{s\prime}f
\bigl(h-H_s'\bigr),
\]
where primes indicate that the future infimum and the height process
are taken w.r.t. the process $X'$ with law $N(\cdot\mid\sup H >h)$,
defined as
\[
X'_t=X_{{\rho_h}+t}-I_{\sigma_h},\qquad 0\le t\le
\tau_h-\rho_h,
\]
where $\rho_h$ is the unique time $s\le\sigma_h$ when $X_s=I_{\sigma
_h}$ and $\tau_h$ is the first time $t\ge\sigma_h$ when
$X_t=I_{\sigma_h}$.

We end the proof with (iii). Thanks to Proposition
\ref{proprwgreat-aunts}, we know that $\langle\rho^{(p)},
f_{p,\gamma _p}\rangle$ has the same law as
\[
\sum_{j=0}^{[\sigma^p_{\gamma_p h}]} \Delta
I^{(p) j}_{\sigma
^p_{\gamma
_p h}}f_{p,\gamma_p}\bigl(\gamma_p
h-H^{(p)}_j\bigr)
\]
conditional on $\max_{1\le j \le T^{(p)}_{-1}} H^{(p)}_j\ge\gamma_p
h$, where $H^{(p)}$ is the height process associated with $W^{(p)}$,
$T^{(p)}_{-1}$ is the first hitting time of $-1$ by $W^{(p)}$, $\sigma
^p_{\gamma_p h}$ denotes the first hitting time of $\gamma_p h$ by
$H^{(p)}$ and $I^{(p)}$ denotes the future infimum process of $W^{(p)}$
stopped at time $\sigma^p_{\gamma_p h}$. If we can prove convergence of
$(p^{-1}W^{(p)}_{[p\gamma_p\cdot\land\sigma^p_{\gamma_p h}]},
\gamma_p^{-1}H^{(p)}_{[p\gamma_p\cdot\land\sigma^p_{\gamma_p h}]})$ under
this conditional law to $(X_{\cdot\land\sigma_h}, H_{\cdot\land
\sigma_h})$ under the measure $N(\cdot\mid\sup H \ge h)$, then by the
generalized continuous mapping theorem (e.g.,~\cite{Kal}, Theorem 4.27)
we will get the convergence of
\[
\sum_{j=0}^{[\sigma^p_{\gamma_p h}]} \frac{\Delta I^{(p) j}_{\sigma
^p_{\gamma_p h}}}{p}f_p
\biggl(h-\frac{H^{(p)}_j}{\gamma_p} \biggr),
\]
which has the law of $p^{-1}\langle\rho^{(p)}, f_{p,\gamma_p}\rangle$,
to $\int_{[0,{\sigma_h}]}d_sI_{\sigma_h}^sf(h-H_s)$ under $N(\cdot\mid
\sup H
>h)$ which has the law of $\langle\rho^0,f\rangle$ thanks to (i).

For (sub)critical $\xi_p$ the convergence of this pair in distribution
on the Skorokhod space of cadlag real-valued paths is a direct
consequence of the results Corollary 2.5.1 and Proposition 2.5.2
already shown in~\cite{DLG}, Section~2.5. To verify that the same holds
for supercritical $\xi_p$ as well, we note that the assumption of
(sub)criticality (hypothesis (H2) in the notation of~\cite{DLG}) is not
crucial in any of the steps of their proof, since the obtained
convergence relies essentially only on the assumption that
$p^{-1}W^{(p)}_{[p\gamma_p\cdot]}$ converges to $X$ (see the comment in
the proof of Theorem 2.2.1 in~\cite{DLG} that any use of subcriticality
in their proof can be replaced by using weak convergence of random
walks). Of course, in the supercritical case, the height process $H$
may drift off to infinity corresponding to the event that the CSB
process $Z$ survives forever, in which case it will code only an
incomplete part of the genealogy of the first lineage which survives
forever. However, for our purposes we only need to consider the
genealogy of the first lineage that survives until time $h$, so this
will not be an impediment for our considerations.

For supercritical $\xi_p$ we still have the convergence of the pair
%
\begin{equation}
\label{eqnconvofWH} \bigl(p^{-1}W^{(p)}_{[p\gamma_p\cdot]},
\gamma_p^{-1}H^{(p)}_{[p\gamma
_p\cdot
]}\bigr) \mathop{\longrightarrow}^{d}_{p\to\infty} (X_\cdot,H_\cdot)
\end{equation}
in distribution on the Skorokhod space of cadlag real-valued paths as
in Theorem~2.3.1 and the first part of Corollary 2.5.1 in~\cite{DLG},
Section 2.5. We now follow the same reasoning as in~\cite{DLG},
Proposition 2.5.2. Let $G^p_{\gamma_p h}= \sup\{s\le\sigma^p_{\gamma_p
h}\dvtx  H^{(p)}_s=0\}$. If we think of $H^{(p)}$ as the height process for
a sequence of independent BGW trees with offspring distribution $\xi
_p$, then $G^p_{\gamma_p h}$ is the initial point of the first BGW tree
in this sequence which reaches a height $\gamma_p h$. Let $(\tilde
{W}{}^{(p)},\tilde{H}{}^{(p)})$ denote the process obtained by conditioning
$(W^{(p)}, H^{(p)})$ on $\max_{1\le j \le T^{(p)}_{-1}} H^{(p)}_j\ge
\gamma_p h$. Then,
\begin{eqnarray*}
&&
\bigl(p^{-1}\tilde{W}{}^{(p)}_{[p\gamma_p\cdot\land\sigma^p_{\gamma_p h}]},
\gamma_p^{-1}\tilde{H}{}^{(p)}_{[p\gamma_p\cdot\land\sigma
^p_{\gamma_p
h}]}\bigr)\\
&&\qquad
\stackrel{d}{=} \bigl(p^{-1}W^{(p)}_{[p\gamma_p(G^p_{\gamma_p
h}+\cdot
)\land\sigma^p_{\gamma_p h}]},
\gamma_p^{-1}H^{(p)}_{[p\gamma
_p(G^p_{\gamma_p h}+\cdot) \land\sigma^p_{\gamma_p h}]}\bigr).
\end{eqnarray*}
Let $G_h=\sup\{s\leq\sigma_h\dvtx H_s=0\}$, and let $(\tilde{X},\tilde{H})$
denote the process obtained by conditioning $(X,H)$ on $\sup H\ge h$, then
\[
(\tilde{X}_{\cdot\land\sigma_h}, \tilde{H}_{\cdot\land\sigma_h}) \stackrel{d}{=}
(X_{(G_h+\cdot)\land\sigma_h},H_{(G_h+\cdot)\land
\sigma_h}).
\]
If we use the Skorokhod representation theorem to assume that the
convergence (\ref{eqnconvofWH}) holds a.s., the same arguments as in
proof Proposition 2.5.2 of~\cite{DLG} imply that $\sigma^{p}_{\gamma_p
h}$ converges a.s. to $\sigma_h$ and $G^p_{\gamma_p h}$ converges a.s.
to $G^p_h$, and hence, $(p^{-1}W^{(p)}_{[p\gamma_p(G^p_{\gamma_p
h}+\cdot)\land\sigma^p_{\gamma_p h}]}, \gamma_p^{-1}H^{(p)}_{[p\gamma
_p(G^p_{\gamma_p h}+\cdot) \land\sigma^p_{\gamma_p h}]})$ converges
a.s. to $(X_{(G_h+\cdot)\land\sigma_h},H_{(G_h+\cdot)\land\sigma
_h})$, proving that
\[
\bigl(p^{-1}\tilde{W}{}^{(p)}_{[p\gamma_p\cdot]},
\gamma_p^{-1}\tilde{H}{}^{(p)}_{[p\gamma_p\cdot]}\bigr)
\mathop{\longrightarrow}^{d}_{p\to
\infty} (\tilde{X}_\cdot,
\tilde{H}_\cdot)
\]
in distribution on the Skorokhod space of cadlag real-valued paths.
From this the desired convergence in distribution of $p^{-1}<\rho
^{(p)}, f_{p,\gamma_p}> $ conditional on $\max_{1\le j \le
T^{(p)}_{-1}} H^{(p)} \ge\gamma_p h$ to $\langle\rho^0,f\rangle$ under
the measure
$N(\cdot\mid \sup H \ge h)$ follows.
\end{pf}

\section{The coalescent point process in the continuous case}

\subsection{Definition of the genealogy}

We now define the analogue of a coalescent point process with
multiplicities for the genealogy (other than the immediately recent) of
an arbitrarily large standing population of a CSB process $Z$. We do so
by first constructing the height process, $H^\star$, for a planar
embedding of CSB trees of arbitrary size descending from an arbitrarily
old ancestor using the continuous version of the great-aunt measure,
$\rho^0$. From $H^\star$ we define coalescence times of two masses from
the standing population in the usual way, that is, from the maximal
depths of the trajectory of $H^\star$.

We construct the height $H_t^\star$ of the individual visited at time
$t$ from: the height $H_t$ defined in (\ref{eqndefheightprocess})
representing the height at which that individual occurs in an excursion
of $X$ above its infimum; plus the height at which this excursion
branches off the unexplored part of the tree. More specifically, let
\[
Y^0_x:=\rho^0\bigl([0,x]\bigr),\qquad x\ge0,
\]
so that
\[
Y^0_x =\beta x + \sum_{t\le x}
\Delta_t,
\]
where\vspace*{1pt} $(t,\Delta_t)$ is a Poisson point measure with
intensity given by $dt\,\pi^{(t)}(dr)$ from Definition~\ref{dfnrho}.
Next, let $L^0$ denote the right-inverse of $Y^0$,
\[
L^0(t):=\inf\bigl\{a\dvtx  Y_a^0>t\bigr\},\qquad t
\ge0.
\]
Then define
\[
H_t^\star:= H_t - L^0(-I_t),
\]
where we remember that $I_t= \inf_{0\le s\le t}X_s$.

This gives a spine decomposition of the genealogy of the
continuous-state branching process associated with $H^\star$ in the
following sense. For the individual visited by the traversal process at
time $t$, the level (measured back into the past, with the present
having level $0$) on the infinite spine at which the subtree containing
this individual branches off is $- L^0(-I_t)$, and the relative height
of this individual within this subtree is $H_t$.

As in the discrete case, we want to display the law of the coalescence
time between successive individuals at generation 0. In this setting,
this corresponds to the maximum depth below 0 of the height process,
between successive visits of 0. Actually, any point at height 0 is a
point of accumulation of other points at height~0, so we discretize the
population at height 0 as follows. We consider all points at height 0
such that the height process between any two of them goes below~$-\vareps$, for some fixed $\vareps>0$, namely, we set $T_0:=0$, and
for any $i\ge1$,
\[
S_i:=\inf\bigl\{t\ge T_{i-1}\dvtx  H_t^\star=-
\vareps\bigr\} \quad\mbox{and}\quad T_i:=\inf\bigl\{t\ge S_{i}\dvtx
H_t^\star=0\bigr\}.
\]
Then the coalescence times $A_1^\vareps, A_2^\vareps,\ldots$ of the
$\vareps$-discretized population are defined as
%
\begin{equation}
\label{eqcoaltimes} A_i^\vareps:=-\inf\bigl
\{H_t^\star\dvtx T_{i-1}\le t\le T_i\bigr
\}.
\end{equation}
As in the discrete case, one of the main difficulties lies in the fact that
the same value of $A_i^\vareps$ can be repeated several times. So we define
%
\begin{equation}
\label{eqmultiplicities} N_i^\vareps:=\#\bigl\{j\ge
i\dvtx A_j^\vareps=A_i^\vareps\bigr\}.
\end{equation}

\subsection{The supercritical case}
\label{subsecsuper}

Actually, the previous definition of genealogy only holds for the
subcritical and critical cases, and a modification needs to be made in
the supercritical case due to possible appearances of branches with
infinite survival times into the future. What we need to do in the
supercritical case is construct a height process $H^\star$ that
corresponds to a CSB tree whose infinite branches have been truncated,
much as in the last chapter of~\cite{L1}. One would be naturally led to
consider the height process of a tree truncated at some finite height.
However, the distribution of such an object is far more complicated
than to truncate (actually, reflect) the associated L\'{e}vy process at
some finite level. On the other hand, the genealogy coded by a L\'{e}vy
process reflected below some fixed level is not easy to read since it
will have incomplete subtrees at different heights. To comply with this
difficulty, we will construct a consistent family of L\'{e}vy processes
$X^\kappa$ reflected below level $\kappa$, so that if $\kappa'>\kappa$,
$X^\kappa$ can be obtained from $X^{\kappa'}$ by excising the subpaths
above $\kappa$. The genealogy coded by the projective limit of this
family is the supercritical L\'{e}vy tree, as is shown in~\cite{L1} in
the case with finite variation. We first give the details of this
construction in the discrete case and then sketch its definition in the
CSB case.
Let $\cal T$ be a planar, discrete, possibly infinite, tree in the
sense that it can have infinite height but all finite breadths. Then
$\cal T$ still admits a Ulam--Harris--Neveu labeling and, as in
Section~\ref{subsecrwr}, we can define $\hat{W}(\varnothing)=0$ and
\[
\hat{W}(v):=\sum_{u\prec v} r(u),\qquad u\not= \varnothing,
\]
where $r(u)$ denotes the number of younger sisters of $u$. Now for any
$\kappa\in\NN$, let ${\cal T}^\kappa$ denote the graph obtained from
$\cal T$ by deleting all vertices $v$ such that $\hat{W}(v)>\kappa$. It
is then easy to see that ${\cal T}^\kappa$ is still a tree, and that if
$v_n^\kappa$ denotes the $n$th vertex of ${\cal T}^\kappa$ in the
lexicographical order and $W_n^\kappa:= \hat{W}(v_n^\kappa)$, then for
any $\kappa' >\kappa$, the path of $W^\kappa$ can be obtained from the
path of $W^{\kappa'}$ by excising the subpaths above $\kappa$. In
other words,
\[
W^\kappa= C_\kappa\bigl(W^{\kappa'}\bigr),
\]
where for any path $\eps$, the functional $C_\kappa$ erasing subpaths
above $\kappa$ can be defined as follows. Let $A^\kappa$ denote the
additive functional
\[
A_n^\kappa(\eps):= \sum_{k=0}^n
\indic{\eps_k\le\kappa},\qquad n\in\NN,
\]
and $a^\kappa$ its right inverse
\[
a_k^\kappa(\eps):= \min\bigl\{n\dvtx  A_n^\kappa(
\eps) > k\bigr\},\qquad k\in\NN.
\]
Then $C_\kappa(\eps):=\eps\circ a^\kappa(\eps)$.

Further, it can be shown that $W^\kappa$ has the law of the random walk
with steps distributed as $\xi-1$, reflected below $\kappa$ and killed
upon hitting $-1$. Since these two properties (consistency by
truncation and marginal distribution) characterize the family of
reflected, killed random walks $(W^\kappa;\kappa\in\NN)$, its
continuous analogue has the following natural definition.

For any $\kappa>0$ let $X^\kappa$ denote a L\'{e}vy process without
negative jumps and Laplace exponent $\psi$ reflected below $\kappa$ and
killed upon reaching 0. The reflection can be properly defined as
follows. Start with the path of a L\'{e}vy process $X$ (with the same
Laplace exponent), set $S_t:=\sup_{0\le s \le t} X_s$ and define
$X^\kappa_t:= X_t$ if at time $t$, $X$ has not yet hit $(\kappa,+\infty
)$, and $X_t^\kappa:=\kappa+X_t -S_t$ otherwise (and kill $X^\kappa$
when it hits 0). The same definition could be done by concatenating
excursions of $X$ away from $(\kappa,+\infty)$ thanks to It\^{o}'s
synthesis theorem.

Now, in the continuous setting, a similar definition of $C_\kappa$ can
be done by adapting the additive functional $A^\kappa$ into
\[
A_t^\kappa(\eps):=\int_0^t
ds \, \indic{\eps_s\le\kappa},\qquad t\ge0.
\]
It is easily seen that for any $\kappa'>\kappa$,
\[
X^\kappa\stackrel{d} {=}C_\kappa\bigl(X^{\kappa'}\bigr).
\]
Kolmogorov's extension theorem then ensures the existence of a family
of processes $(X^\kappa;\kappa>0)$ defined on the same probability
space, satisfying pathwise the previous equality, and with the right
marginal distributions (reflected, killed L\'{e}vy processes).

All results stated in the remainder of the paper hold in the
supercritical case if we replace the L\'{e}vy process $X$ by the
projective limit of $X^\kappa$ as $\kappa\to\infty$. More simply, we
can construct the same consistent family $(e^\kappa;\kappa>0)$ for the
excursion of $X-I$ away from 0. Then it is sufficient to replace each
excursion of $X-I$ drifting to $\infty$ (there are finitely many of
them on any compact interval of local time) by a copy of $e^\kappa$ for
some sufficiently large $\kappa$. More precisely, for the excursion
corresponding to an infimum equal to $-x$ ($x$ is the index of the
excursion in the local time scale), we need that the modified excursion
$e^\kappa$ be such that all heights below $h:=L^0(x)$ be visited. In
other words, one has to choose $\kappa$ large enough so that for any
$\kappa'>\kappa$, the occupation measure of the height process of
$e^{\kappa'}$ restricted to $[0,h]$ remains equal to that of $e^\kappa$.

In order not to overload the reader with technicalities that are away
from the core question of this paper, we chose not to develop this
point further, and we leave it to the reader to modify the proofs of
the next subsection in the obvious way at points where the
supercritical case has to be distinguished.


%

\subsection{Law of the coalescent point process}
Now that we have the height process for the infinite CSB tree with an
arbitrarily old genealogy and the ancestral coalescence times with
multiplicities describing genealogy of its standing population, we
proceed to describe their law. Furthermore, as was our main goal, we
define an analogue of the coalescent point process with multiplicities
$(B_i;i\ge0)$ for the CSB tree.

Our results show similarities with the results of Duquesne and Le Gall
in~\cite{DLG}, Section 2.7, for the law of the reduced tree under the
measure $N(\cdot\mid\sup H \ge T)$. However, our presentation is
quite different for a number of reasons. First, we do not characterize
the branching times in terms of the Markov kernel of the underlying CSB
process, but rather treat these times as a sequence. Second, the
branching tree is allowed to be supercritical. Third, the tree may have
an infinite past.

Recall the definition of coalescence times of individuals in the
discretized standing population $(A^\vareps_i;i\ge1)$ and their
multiplicites $(N^\vareps_i;i\ge1)$ from (\ref{eqcoaltimes}) and
(\ref{eqmultiplicities}).
%
\begin{theorem}
\label{thmlawofAeps1}
The joint law of $(A_1^\vareps, N_1^\vareps)$ is given by
%
\begin{equation}
\label{eqntailofAeps} \PP\bigl(A_1^\vareps> x\bigr)=
\frac{\tilde{\psi}(v(x))}{\tilde{\psi
}(v(\vareps
))},\qquad x\ge\vareps,
\end{equation}
where $\tilde{\psi}(u):=\psi(u)/u$ for any $u>\eta$. In addition, for
any $n\ge2$ and $x\ge\vareps$,
\[
\PP\bigl(A_1^\vareps\in dx, N_1^\vareps=n
\bigr)/dx= \PP\bigl(A_1^\vareps>x\bigr)v(x)^n\int
_{(0,\infty)}\Lambda(dz) e^{-v(x)z}\frac{z^{n+1}}{(n+1)!},
\]
whereas for any $x\ge\vareps$,
\[
\PP\bigl(A_1^\vareps\in dx, N_1^\vareps=1
\bigr)/dx= \PP\bigl(A_1^\vareps>x\bigr)v(x) \biggl(\beta+\int
_{(0,\infty)}\Lambda(dz) e^{-v(x)z}\frac{z^2}{2} \biggr).
\]
\end{theorem}
%
\begin{rem}
The formula giving the joint law of $(A_1^\vareps,N_1^\vareps)$ can
also be expressed (see the calculations in the proof)
as follows:
\begin{eqnarray*}
&&\PP\bigl(A_1^\vareps\in dx, N_1^\vareps=n
\bigr)/dx\\
&&\qquad= \PP\bigl(A_1^\vareps>x\bigr) \biggl(\beta
\indic{n=1}+\int_{(0,+\infty)} \pi^{(x)}(dr) e^{-rv(x)}
\frac
{(rv(x))^n}{n!} \biggr),
\end{eqnarray*}
showing that when $\beta=0$, $N_1^\vareps$ actually follows a mixed
Poisson distribution.
\end{rem}

\begin{rem} Similarly, as in the discrete case, observe that for
subcritical trees $[\eta=0$ and $\psi'(0+)>0]$, coalescence times can
take the value $+\infty$, corresponding to the delimitation of
quasi-stationary subpopulations with different (infinite) ancestor
lineages. Indeed, from the previous statement,
\[
\PP\bigl(A_1^\vareps=+\infty\bigr)= \frac{\psi'(0+)}{\tilde{\psi
}(v(\vareps))}.
\]

In addition, the event $\{A_1^\vareps=+\infty\}$ is the event that the
whole (quasi-stationary) population has coalesced by time $\vareps
$. In other words, if $V$ denotes the coalescence time of a
quasi-stationary population, also referred to as the time to most
recent common ancestor, then
\[
\PP(V>x)= \frac{\psi'(0+)}{\tilde{\psi}(v(x))}.
\]
In the discrete case, the study of $V$ will be done in Section
\ref{subsecdisintegration}, in a slightly more detailed statement
Proposition~\ref{propmrcaqsd}.
\end{rem}
\begin{pf*}{Proof of Theorem \protect\ref{thmlawofAeps1}}
First notice that $S_1$ is also the first hitting time of $-Y^0_\vareps
$ by $X$. Denote by $(s,e_s)$ the excursion process of $X-I$\vadjust{\goodbreak} away from
0, where $-I$ serves as local time. Then the event $\{A_1^\vareps>x\}$
is the event that the excursions $(e_s; Y^0_\vareps\le s \le Y^0_x)$
all have $\sup H(e_s)-L^0_s<0$. As a consequence,
\[
\PP\bigl(A_1^\vareps>x\bigr)=\EE\exp-\sum
_{Y^0_\vareps\le s \le Y^0_x}\chi_{\{
\sup H(e_s)<L^0_s\}},
\]
where $\chi_{\{A\}}=0$ on $A$ and $\chi_{\{A\}}=+\infty$ on $^c A$, so
that using the exponential formula for the excursion point process,
\[
\PP\bigl(A_1^\vareps>x\bigr)=\EE\exp-\int
_{[Y^0_\vareps, Y^0_x]} ds\, N\bigl(\sup H>L^0_s\bigr)=
\EE\exp-\int_{[Y^0_\vareps, Y^0_x]} ds\, v\bigl(L^0_s
\bigr).
\]
Now, since $L^0$ is the right-inverse of $Y^0$, and because of the
definition of $Y^0$ (Definition~\ref{dfnrho}) in terms of the
Lebesgue measure and the Poisson point process $(u,\Delta_u)$, we get
\[
\int_{[Y^0_\vareps, Y^0_x]} ds\, v\bigl(L^0_s\bigr) =
\beta\int_{[\vareps,
x]}du\, v(u)+\sum_{\vareps\le u\le x}v(u)
\Delta_u,
\]
so that
%
\begin{eqnarray}
\label{eqAepsfromPP}
&&
\PP\bigl(A_1^\vareps>x\bigr)\nonumber\\[-8pt]\\[-8pt]
&&\qquad=\exp-
\biggl(\beta\int_{[\vareps, x]}du\, v(u)+\int_{[\vareps, x]}du
\int_{[0,\infty)}\pi^{(u)}(dr) \bigl(1-e^{-rv(u)}
\bigr) \biggr).\nonumber
\end{eqnarray}
We compute the second term inside the exponential thanks to the
Fubini--Tonelli theorem as
\begin{eqnarray*}
\int_{[0,\infty)}\pi^{(u)}(dr) \bigl(1-e^{-rv(u)}
\bigr)&=&\int_{[0,\infty)}dr\, e^{v(u)r} \int
_{(r,\infty)} e^{-v(u)z} \Lambda(dz) \bigl(1-e^{-rv(u)}
\bigr)
\\
&=&\int_{(0,\infty)} \Lambda(dz) e^{-v(u)z}\int
_{(0,z)}dr\, \bigl(e^{v(u)r}-1 \bigr)
\\
&=&\frac{1}{v(u)}\int_{(0,\infty)} \Lambda(dz) e^{-v(u)z}
\bigl(e^{v(u)z}-1-v(u)z \bigr).
\end{eqnarray*}
Finally, we get
\[
\PP\bigl(A_1^\vareps>x\bigr)=\exp-\int_{[\vareps, x]}du
\,F\circ v(u),
\]
where
\[
F(\lbd):=\beta\lbd+ \frac{1}{\lbd}\int_{(0,\infty)} \Lambda(dz)
\bigl(1-e^{-\lbd z}-\lbd ze^{-\lbd z} \bigr),\qquad \lbd>\eta.
\]
But elementary calculus shows that
\[
F(\lbd) = \psi'(\lbd) - \tilde{\psi}(\lbd),\qquad \lbd>\eta,
\]
so that by the change of variable $y=v(u)$, $u=\phi(y)$, we get
\begin{eqnarray*}
\PP\bigl(A_1^\vareps>x\bigr)&=&\exp-\int
_{[v(x), v(\vareps)]}\frac{dy}{\psi(y)}\, F(y)
\\
&=&\exp-\int_{[v(x), v(\vareps)]}dy \,\biggl(\frac{\psi'(y)}{\psi(y)}-
\frac{1}{y} \biggr)
\\
&=&\frac{\tilde{\psi}(v(x))}{\tilde{\psi}(v(\vareps))},
\end{eqnarray*}
which shows the first part of the theorem.

As for the second part, the event $\{A_1^\vareps\in dx, N_1^\vareps=
n\}$ is the event that the excursions $(e_s; Y^0_\vareps\le s \le
Y^0_{x-dx})$ all have $\sup H(e_s)-L^0_s<0$, and $n$ is the number of
excursions $(e_s; Y^0_{x-dx}\le s\le Y^0_{x})$ (for which $L_s^0\in
dx$) such that $\sup H(e_s) \ge x$. Next, notice that
$Y^0_{x}-Y^0_{x-dx} =\beta \,dx+ (Y^0_{x}-Y^0_{x-})\indic{E^0(dx)}$,
where $E^0(dx)$ is the event that $Y^0$ has a jump in the interval
$(x-dx,x)$. As a consequence, by the compensation formula applied to
the Poisson point process of jumps of $Y^0$,
\[
\PP\bigl(A_1^\vareps\in dx, Y^0_{x}-
Y^0_{x-}\in dr, N_1^\vareps=n\bigr)=
\PP\bigl(A_1^\vareps>x\bigr)\,dx\,\pi^{(x)}(dr)
e^{-rv(x)}\frac{(rv(x))^n}{n!},
\]
since $N(\sup H\ge x)=v(x)$. Also
\[
\PP\bigl(A_1^\vareps\in dx, Y^0_{x}=
Y^0_{x-}, N_1^\vareps=n\bigr)=0,
\]
if $n\ge2$, whereas
\[
\PP\bigl(A_1^\vareps\in dx, Y^0_{x}=
Y^0_{x-}, N_1^\vareps=1\bigr)= \PP
\bigl(A_1^\vareps>x\bigr)\beta \,dx \, v(x).
\]
As a consequence, for any $n\ge2$,
\begin{eqnarray*}
\PP\bigl(A_1^\vareps\in dx, N_1^\vareps=n
\bigr)&=& dx\,\PP\bigl(A_1^\vareps>x\bigr)\frac
{v(x)^n}{n!}
\int_{[0,\infty)}\pi^{(x)}(dr) e^{-rv(x)}r^n
\\
&=& dx\,\PP\bigl(A_1^\vareps>x\bigr)\frac{v(x)^n}{n!}\int
_{[0,\infty)}dr\, r^n\int_{(r,\infty)}
\Lambda(dz) e^{-zv(x)}
\\
&=& dx\,\PP\bigl(A_1^\vareps>x\bigr) v(x)^n\int
_{(0,\infty)}\Lambda(dz) e^{-v(x)z}\frac{z^{n+1}}{(n+1)!},
\end{eqnarray*}
which is the desired result. The last result can be obtained by the
same calculation. Summing over $n$ yields
\begin{eqnarray*}
\PP\bigl(A_1^\vareps\in dx\bigr)&=&dx\,\PP
\bigl(A_1^\vareps>x\bigr) \\
&&{}\times\biggl(\beta v(x)+\frac
{1}{v(x)}
\int_{(0,\infty)}\Lambda(dz) \bigl(1-e^{-v(x)z}-v(x)ze^{-v(x)z}
\bigr) \biggr)
\\
&=&dx\,\PP\bigl(A_1^\vareps>x\bigr) F\circ v(x)
\end{eqnarray*}
as expected, since $\PP(A_1^\vareps>x)=\exp-\int_{[\vareps, x]}du\,
F\circ v(u)$.\vadjust{\goodbreak}
\end{pf*}

Finally, we define the coalescent point process with multiplicities for
the $\vareps$-discretized standing population. At the end of
Section~\ref{subsecdiscrcoal} we gave an alternative
characterization of $(B_k;k\ge0)$ for the BGW tree as a sum of point
masses [see (\ref{eqaltcoalesc})] whose multiplicities record the
number of their future appearances as coalescent times, until the first
future appearance of a larger coalescence time. Since in the CSB tree
we do not have a process analogous to $(D_i;i\ge1)$, we actually
define the $\vareps$-discretized coalescent point process with
multiplicities from this angle.

For each $i\ge1$ and $k\ge i$ we define $N_{ik}^\vareps$ as the
residual multiplicity of $A_{i}^\vareps$ at time~$k$,
\[
N_{ik}^\vareps:=\#\bigl\{j\ge k\dvtx A_j^\vareps=A_i^\vareps
\bigr\},
\]
so $N_{ii}^\vareps=N_i^\vareps$. Next define for each $k\ge1$ the
random finite point measure $B_k^\vareps$ on $(\vareps,+\infty]$ as
\[
B_k^\vareps:=\sum_{i\le k}N_{ik}^\vareps
\delta_{A_i^\vareps
}\mathbf{1}_{\{A^\vareps_i< A^\vareps_{i'},\forall i'<i\dvtx N^\vareps
_{i'k}\neq
0\}},
\]
where $\delta$ denotes a Dirac measure, and by convention let
$B_0^\vareps$ be the null measure.

Recall that $\supp$ is a mapping from a point measure on $\RR$ to the
minimum of its support. The following result provides the law for the
coalescent point process with multiplicities of the $\vareps
$-discretized population in the CSB tree.

\begin{theorem}
\label{thmlawofcontscoal}
The sequence of finite point measures $(B_i^\vareps;i\ge0)$ is a
Markov chain started at the null measure. For any finite point measure
$b=\sum_{j\ge1}n_j\delta_{a_j}$, with $n_j\in\mathbb{N}$ and
$a_j\in
\mathbb{R^+}$, 
the law of $B_{k+1}^\vareps$ conditionally given $B_k^\vareps=b$, is
given by the following transition.
Let $a_1:=\supp(b), b^*:=b-\delta_{a_1}, a_1^*:=\supp(b^*)$.
Let $(A,N)$ be a r.v. with values in $(\vareps,+\infty]\times\NN$
distributed as $(A_1^\vareps, N_1^\vareps)$. Then
\[
B^\vareps_{k+1}:= \cases{ b^*+N\delta_A, &\quad if
$A< a_1^*$,
\cr
b^*, &\quad otherwise. }
\]
%
In addition, as in the discrete case, the coalescent point process can
be deduced from $B^\vareps$ as
\[
A_k^\vareps=\supp\bigl(B^\vareps_k\bigr).
\]
\end{theorem}
%
\begin{rem}
When $X$ is a diffusion, the measure $\rho^0$ has no atoms, so there
are no repeats of coalescence times ($N=1$ a.s.). In this case, at
each step of the chain there is only one nonzero mass, that is,
$\forall k\ge1, B_k^\vareps=\delta_{A_k^\vareps}$. The
previous\vspace*{1pt}
statement shows that the sequence starts with $B_1^\vareps=\delta
_{A_1^\vareps}$ where $A^\vareps_1$ is distributed as $A$, and in every
transition the single point mass from the previous step is erased and a
new point mass $\delta_{A_{k+1}^\vareps}$ is added with independently
chosen $A_{k+1}^\vareps$ distributed as $A$. So the random variables
$(A_i^\vareps)_{i\ge1}$ are i.i.d.\vadjust{\goodbreak}

On the other hand, when $X$ is a diffusion (and only in that case;
see~\cite{DLG}), the height process $H$ is a Markov process.
The coalescence times $(A^\vareps_i;i\ge1)$ are just the depths of
excursions of $H$ (with depth greater than $\varepsilon$) below some
fixed level. This again implies that the $(A_i^\vareps)_{i\ge1}$ are i.i.d.
Moreover, by taking $\vareps\to0$ in equation (\ref{eqntailofAeps}),
one can compute the intensity measure $\mu$ of the Poisson point
process of excursion depths [its tail is given by $\bar{\mu
}(x)=\tilde
{\psi}(v(x))$].
In particular, in the Brownian case, it is known that the height
process is (reflected) Brownian, so that $\mu(dx)$ is proportional
to $x^{-2}\,dx$ (see~\cite{AP,P}).

For the BGW tree, the discrete analogue of the height process $H$
is, again, in general not a Markov process, except in special cases of
the offspring distribution, namely, the only exceptions are when $\xi$
is linear fractional (see Section~\ref{subseclinfrac} for
definition). In these cases we will also observe that the coalescence
times \mbox{$(A_i;i\ge1)$} are i.i.d. (see Proposition~\ref{propcpplf}).
\end{rem}
%
\begin{rem}
Even in general, that is, in the presence of multiplicities, we know
that there exists a coalescent point process whose truncation at level
$\vareps$ is the sequence $(A^\vareps_i;i\ge1)$. It is the process of
excursion depths of the height process in the local time scale.
However, the question of characterizing the distribution of this
coalescent point process (without truncation) remains open. A
natural idea would be to use a Poisson point process with intensity
measure $dt \,\sum_{n\ge1}\nu(dx, n) \delta_n$, where
\[
\nu\bigl(dx, \{n\}\bigr)/dx= \tilde{\psi}\bigl(v(x)\bigr) v(x)^n
\biggl(\beta\indic{n=1}+\int_{(0,\infty)}\Lambda(dz) e^{-v(x)z}
\frac
{z^{n+1}}{(n+1)!} \biggr).
\]
Indeed, since
$(A_1^\vareps, N_1^\vareps)$ has the same law as the first atom of
this point process conditional on its first component being larger than
$\vareps$, it is easy to embed the sequence $(A_i^\vareps, i\ge1)$
into the atoms of this Poisson point process by truncating atoms larger
than those whose multiplicity has not yet been exhausted. However,
preliminary calculations indicate that convergence of this new point
process as $\vareps\to0$ is not evident.
\end{rem}
\begin{pf*}{Proof of Theorem \protect\ref{thmlawofcontscoal}}
We need to introduce some notation. For any stochastic process $W$
admitting a height process in the sense of (\ref{eqndefheightprocess}),
we denote by $H^W$ its height process, that is,
\[
H^W_t:= \lim_{k\to\infty}\frac{1}{\vareps_k}\int
_0^t \indic{W_s<I_{W,t}^s+
\vareps_k} \,ds,
\]
where
\[
I_{W,t}^s:= \inf_{s\le r\le t} W_r.
\]
Further, if $W$ has finite lifetime, denoted $T$, we denote by $\rho^W$
the great-aunt measure associated with $W$, in the sense that
\[
\bigl\langle\rho^W,f\bigr\rangle:=\int_{[0,T]}d_sI_{W,T}^s
f\bigl(H_T^W-H_s^W\bigr).\vadjust{\goodbreak}
\]
In particular, if $W$ is the L\'{e}vy process killed at $\sigma_h$
under $N(\cdot\mid\sup H> h )$, then we know from Theorem \ref
{thmpropertiesrho} that $\rho^W$ and the trace of $\rho^0$ on $[0,h]$ are
equally distributed. Now if $\mu$ is a positive measure on $\RR_+$ and
$L^\mu$ denotes the right-inverse of the nondecreasing function
$x\mapsto\mu([0,x])$, and if $I^W$ denotes the past infimum of the
shifted path $W-W_0$, we denote by $\Phi(\mu, W)$ the generalized
height process
\[
\Phi_t(\mu, W):=H^W_t- L^\mu
\bigl(-I^W_t\bigr).
\]
In particular, in the (sub)critical case, $H^\star$ is distributed as
$\Phi(\rho^0, X)$. Finally, recall that $\kil$ denotes the killing
operator and $\theta$ the shift operator, in the sense that $X\circ
\theta_t = (X_{t+s};s\ge0)$ and $X\circ\kil_t$ is $(X_{s \land
t};s\ge0)$. For any path $(X_s;s\ge0)$ and any positive real number
$t$, if $W:=X\circ\kil_t$ and $W'=X\circ\theta_t$, then it is not
hard to see that for any $s\ge0$,
%
\begin{equation}
\label{eqnnothard} H^X_{t+s} =\Phi_s
\bigl(\rho^W, W'\bigr)+ H^W_t.
\end{equation}

In particular, applying this to $X$ under $N(\cdot\mid\sup H> h )$ and
to $t=\sigma_h$, the strong Markov property at $\sigma_h$ yields
%
\begin{equation}
\label{eqnweget} H^X\circ\theta_{\sigma_h} \stackrel{d}
{=} \Phi(\rho, X\circ\theta_{\sigma_h})+ h,
\end{equation}
where $\rho$ is a copy of $\rho^0$ independent of $X$.

Now we work conditionally on $(A_1^\varepsilon, N_1^\varepsilon
)=(a,n)$. Recall that $S_1=\inf\{t\ge0\dvtx  H_t^\star=-\varepsilon\}$ and
$T_1=\inf\{t\ge S_1\dvtx  H_t^\star=0\}$. Set $V_1:=T_1$ and define
recursively for $i=1,\ldots,n$,
\begin{eqnarray*}
U_i&:=&\sup\bigl\{t<V_i\dvtx  H_t^\star=-a
\bigr\},
\\
W_{i}&:=&\inf\bigl\{t>V_i\dvtx  H_t^\star=-a
\bigr\},
\\
V_{i+1}&:=&\inf\bigl\{t>W_i\dvtx  H_t^\star=0
\bigr\}.
\end{eqnarray*}
%
Last, define
\[
W_0:=\inf\bigl\{t>0\dvtx  X_t=-Y^0_{a-}
\bigr\} \quad\mbox{and}\quad U_{n+1}:=\inf\bigl\{t>0\dvtx  X_t=-Y^0_a
\bigr\}.
\]
As seen in the proof of Theorem~\ref{thmlawofAeps1}, the subpaths
$e^i:=(X_{t+U_i}-X_{U_i}; 0\le t \le W_{i}-U_i)$ are the $n$ excursions
of $X$ above its infimum whose height reaches level~$a$, so that for
all $t\in[W_0, U_{n+1}]$, $I_t\in(-Y^0_{a-}, Y^0_a]$ and $H_t^\star=
H_t -a$. An application of the strong Markov property yields the
independence of the $n$ subpaths $e^i$ and the fact that they are all
distributed as $N(\cdot\mid\sup H\ge a)$. Also,
$X':=(X_{t+U_{n+1}}-X_{U_{n+1}}; t\ge0)$ is a copy of $X$ independent
from all the previous subpaths and from~$Y^0$. Notice that for all
$t\in [W_0, U_{n+1}]$, $H_t^\star\ge-a$. As a consequence, by
continuity of the height process, $H^\star$ takes the value $-a$ at all
points of the form $U_i$ and $W_i$, it hits 0 only in intervals of the
form $[V_i, W_i)$ and takes values in $[-a,0)$ on all intervals of the
form $[W_i, U_{i+1}]$. As a consequence, if we excise all the paths of
$H^\star$ on intervals of the form $(W_i, U_{i+1})$, $i=0,\ldots,n$, we
will still get the same coalescent point process. Also notice that
those paths are independent of the $n$ excursions $e^i$ and only depend
on $Y^0$ through $\Delta Y_a^0$. As a consequence, recalling the
notation in Section~\ref{subsecsuper}, if $\mathcal{A}^\vareps$ denotes
the mapping that takes a height process to its sequence of
$\vareps$-discretized coalescent levels, that is,
$\mathcal{A}^\vareps\dvtx H^\star\mapsto(A_1^\vareps ,A_2^\vareps,\ldots)$,
where $A_1^\vareps,A_2^\vareps,\ldots$ are defined by
(\ref{eqcoaltimes}), then we have
\[
\mathcal{A}^\vareps\bigl(H^\star\bigr) = \bigl(a,
\mathcal{A}^\vareps\bigl(H\bigl(e^1\bigr)-a\bigr),a,\ldots,
\mathcal{A}^\vareps\bigl(H\bigl(e^{n-1}\bigr)-a\bigr), a,
\mathcal{A}^\vareps\bigl(H_n'\bigr)\bigr),
\]
where $H_n'$ is the concatenation of $H(e^n)-a$ and of $\Phi(\rho^0\circ
\theta_a, X')$, writing $\rho^0\circ\theta_a$ for the measure
associated with the jump process $(Y^0_{s+a}-Y^0_a;s\ge0)$. First
observe that, as long as the coalescent point process is concerned, we
can again excise the parts of each of the $n$ subpaths in the previous
display going from height $-a$ to height $0$. This amounts to
considering excursions $e^i$ only from the first time $\sigma_a$ they
reach height $a$. But recall from (\ref{eqnweget}) that $H(e)\circ
\theta_{\sigma_a} -a\stackrel{d}{=} \Phi(\rho, e\circ\theta_{\sigma
_a})$, which is distributed as the process $H^\star$ killed upon
reaching $-a$. Second, by the same argument as stated previously,
erasing the part of $H_n'$ before its first hitting time of $0$, we get
the concatenation of a copy of $H^\star$ killed upon reaching $-a$ and
of $\Phi(\rho^0\circ\theta_a, X')$, where we remember that $X'$ is an
independent copy of $X$. The result has the law of $\Phi(\rho^0, X)$,
that is, it has the law of $H^\star$. In conclusion, this gives the
following conditional equality in distribution:
\[
\mathcal{A}^\vareps\bigl(H^\star\bigr) = \bigl(a,
\mathcal{A}^\vareps\bigl(H_1^\star\bigr),a,\ldots,
\mathcal{A}^\vareps\bigl(H_{n-1}^\star\bigr), a,
\mathcal{A}^\vareps\bigl(H_n''
\bigr)\bigr),
\]
where $H_n''$ is a copy of $H^\star$ and the $H_i^\star$ are
independent copies of $H^\star$ killed upon reaching $-a$, all
independent of $H_n''$.

Now observe that since the law of $A_1^\vareps$ is absolutely
continuous, the branch length~$a$ will occur exactly $n$ times [in
particular, it will not appear in $\mathcal{A}^\vareps(H_n'')$]. Also,
because all the heights between successive occurrences of $a$ are
smaller than~$a$, the following conditional equality can be inferred
from the previous display and the definition of the sequence $B^\vareps
:=(B^\vareps_i)$ of point measures
\[
B^\vareps= \bigl(0,n\delta_a, (n-1)\delta_a +
B^1, (n-2)\delta_a,\ldots, \delta_a+B^{n-1},
\delta_a, B'\bigr),
\]
where $B'$ is a copy of $B^\vareps$ and the $B^i$ are independent
copies of the sequence of point measures associated with a coalescent
point process $A^\vareps$ killed at its first value greater than $a$
(in the usual sense that this value is not included in the killed
sequence). In passing, an induction argument on the cardinal of the
support of $B^\vareps$ shows that $A^\vareps=\supp(B^\vareps)$.
Also, each $B^i$ is reduced to a sequence of length 1 with single value
0, with probability $\PP(A^\vareps>a)$. Otherwise, it starts in the
state~$N\delta_A$, where $(A,N)$ has the law of $(A_1^\vareps,
N_1^\vareps)$ conditional on $A_1^\vareps<a$. Note that the sequence
following the state $\delta_a$ [which corresponds to the case when
$\supp(b)\not=\supp(b^*)$], has the same law as $B^\vareps$, dropping
the dependence upon $a$.

Now assume that $B_k^\vareps=\sum_{j\ge1}n_j \delta_{a_j}$. The
following conclusions follow by induction from the last assertion. If
$n_1 \ge2$ (i.e., $a_1^*=a_1$) and $K$ is the first time after $k+1$
that $a_1$ has its multiplicity decreased, then the sequence
$(B_i^\vareps; k+1\le i < K )$ is independent of $(B_i^\vareps;i\le k)$
and has the law of the sequence associated with a coalescent point
process killed at its first value greater than $a_1$. If $n_1=1$ (i.e.,
$a_1^*=a_2$) and $K$ is the first time (after $k$) that $a_2$ has its
multiplicity decreased, then the sequence $(B_i^\vareps; k+1\le i < K
)$ is independent of $(B_i^\vareps;i\le k)$ and has the law of the
sequence associated with a coalescent point process killed at its first
value greater than $a_2$.
In particular, $B_{k+1}^\vareps= B_k^\vareps- \delta_{a_1}$ with
probability $\PP(A^\vareps>a_1^*)$, and $B_{k+1}^\vareps=
B_k^\vareps
- \delta_{a_1}+N\delta_A$, with probability $\PP(A^\vareps< a_1^*)$,
where $(A,N)$ has the law of $(A_1^\vareps, N_1^\vareps)$ conditional
on $A_1^\vareps<a_1^*$.
%
%
%
\end{pf*}

\subsection{Convergence of the coalescent point process}
We now present the theorem that connects the discrete case coalescent
process, based on the offspring distribution $\xi_p$, to the continuous
case coalescent process, based on the associated branching mechanism
$\psi$.
Let us assume, as in Section~\ref{subsecconvofg-a}, that for some
sequence $(\gamma_p, p\ge1), \gamma_p\to\infty$ as $p\to\infty$
and a
sequence of r.v. $(\xi_p,p\ge1)$ such that the rescaled BGW process
started at $[px]$ with offspring distribution $\xi_p$ when re-scaled to
$p^{-1}Z^{(p)}_{[\gamma_p \cdot]}$, converges in Skorokhod space to a
CSB process $Z$ with branching mechanism $\psi$ started at $x$.

In order to obtain convergence for the coalescent processes we need to
define a discretized version of the discrete case process by
considering only the individuals whose coalescent times are greater
than $\gamma_p\vareps$ for some fixed $\vareps>0$. Recall that for a
point measure $b$, $\supp(b)$ denotes the minimum of its support. Start
with the sequence of finite point measures $(B_i;i\ge0)$ whose law is
given by Theorem~\ref{thmcpmp}. Let $\tau_0:=0$, and for any $i\ge
1$, let
\[
\tau_i:=\inf\bigl\{n\ge\tau_{i-1}\dvtx  \supp(B_n)
\ge\gamma_p\vareps\bigr\}.
\]
Define the $\gamma_p\vareps$-discretized process of point measures
\[
\bigl(B^{p,\vareps}_i;i\ge0\bigr):=(B_{\tau_i};i\ge0).
\]
Let $\mathcal{B}= \{\sum_{i=1}^n b_i\delta_{a_i}\dvtx  n\in\NN,
b_i\in\NN, a_i\in\RR_+ \}$ be the space of all finite point mass
measures on
$\RR_+$, equipped with the usual vague topology.
Let ${\cal R}^p\dvtx \mathcal{B}\mapsto\mathcal{B}$ be a function
re-scaling the point mass measures so that
\[
{\cal R}^p\dvtx  \sum_{i=1}^nb_i
\delta_{a_i} \mapsto\sum_{i=1}^nb_i
\delta_{\gamma_p^{-1}{a_i.}}
\]
Since $(B^{p,\vareps}_i;i\ge0)$ is a Markov chain on $\mathcal{B}$, it
is a random element of $\mathcal{B}^{\{0,1,2,\ldots\}}$, equipped with
the product of vague topologies on $\mathcal{B}$.
%
\begin{theorem}
\label{thmconvcpp}
The sequence of re-scaled discretized Markov chains\break $({\cal
R}^p(B^{p,\vareps}_i); i\ge0)$ 
converges in distribution on the space $\mathcal{B}^{\{0,1,2,\ldots\}}$
to the Markov chain $(B^\vareps_i; i\ge0)$ whose law is given by
Theorem~\ref{thmlawofcontscoal} as $p\to\infty$.\vadjust{\goodbreak}
\end{theorem}
\begin{pf}
Note that the initial values for the sequence ${\cal R}^p(B^{p,\vareps
}_0)$, as well as for the limit $B^\vareps_0$, are simply null
measures. In order to describe the transition law of the discretized
process $(B^{p,\vareps}_i; i\ge0)$, condition on its value at step $i$
\[
B_{\tau_i}=b \qquad\mbox{with } a_{1}=\supp(b)\geq
\gamma_p\vareps.%
\]
Condition further on the values of $\tau_i$ and $\tau_{i+1}$ for the
unsampled process $(B_k; k\ge1)$. With $b^*=b-\delta_{a_{1}}$, we have
that for all $1\le j\le\tau_{i+1}-\tau_i-1$,
\[
B_{\tau_i+j}=b^*+\sum_{j\ge1}n_{j}
\delta_{a_{j}} \qquad\mbox{with } a_{j}<\gamma_p\vareps,
\forall j\ge1
\]
since, by definition of $\tau_{i+1}$ for $1\le j\le\tau_{i+1}-\tau
_i-1$, the smallest mass in $B_{\tau_i+j}$ must be smaller than
$\gamma_p\vareps$, and in each step only the weight of the smallest
mass is decreased.

On the event $\tau_{i+1}=\tau_i+1$, the transition rule from
Theorem~\ref{thmcpmp} for step $\tau_i$ to $\tau_{i+1}$ gives that,
for $a_1^*=\supp(b^*)$,
%
\begin{equation}
\label{eqtrsamp} B_{\tau_{i+1}}= \cases{ b^*+N^{p,\vareps}
\delta_{A^{p,\vareps}}, &\quad if $A^{p,\vareps
}<a_{1}^*$ and
$A^{p,\vareps}\neq a_{1}$,
\cr
b^*, &\quad otherwise, }
\end{equation}
where $(A^{p,\vareps},N^{p,\vareps})$ is distributed as $(A_1,\zeta
'_{A_1})$ conditional on $A_1\ge\gamma_p\vareps$. The conditioning in
this law follows from the definition of $\tau_{i+1}$ as the first time
after $\tau_i$ for which $\supp(B^{p,\vareps}_{i+1})\geq\gamma_p\vareps$.

On the\vspace*{1pt} event $\tau_{i+1}>\tau_i+1$, 
since $\supp(B^{p,\vareps}_{i+1})\geq\gamma_p\vareps$, by step
$\tau_{i+1}-1$ all of the masses from $\sum_{j\ge1}n_{j}\delta_{a_{j}}$
must have been eliminated except for one mass $a$ that is smaller than
$\gamma_p\vareps$ whose weight at this step is $1$, so
\[
B_{\tau_{i+1}-1}=b^*+\delta_{a} \qquad\mbox{with } a<\gamma_p
\vareps.
\]
Since the smallest mass in $B_{\tau_{i+1}-1}$ is $\delta_{a}$, the
transition rule from Theorem~\ref{thmcpmp} for step $\tau_{i+1}-1$
to $\tau_{i+1}$ gives that $B^*_{\tau_{i+1}-1}=B_{\tau
_{i+1}-1}-\delta_a=b^*$, so $\supp(B^*_{\tau_{i+1}-1})=\supp
(b^*)=a_1^*$. Also, a new
mass $(A,N)$ is added only if $A<a_1^*$ and $A\neq a=\supp(b^*+\delta_a)$.
Note that we must also have $A\ge\gamma_p\vareps$ as in the next step
$\supp(B_{\tau_{i+1}})\geq\gamma_p\vareps$, so the added mass is again
distributed as $(A_1,\xi_{A_1}')$ conditional on $A_1\ge\gamma_p\vareps$.
Integrating over possible values for $B_{\tau_{i+1}-1}$, $\tau_i$, and
$\tau_{i+1}$, we have that, conditionally given $B^{p,\vareps}_i=b$,
the transition rule for $B^{p,\vareps}_{i+1}$ is
\[
B_{\tau_{i+1}}= \cases{ B_{\tau_i}^*+N^{p,\vareps}
\delta_{A^{p,\vareps}}, &\quad if $A^{p,\vareps}<a_{1}^*$ and
$A^{p,\vareps}\neq a_{1}$,
\cr
B_{\tau_i}^*, &\quad otherwise, }
\]
where $(A^{p,\vareps},N^{p,\vareps})$ is distributed as $(A_1,\zeta
'_{A_1})$ conditional on $A_1\ge\gamma_p\vareps$.

For the re-scaled process ${\cal R}^p(B^{p,\vareps})$, the transition
rule for ${\cal R}^p(B^{p,\vareps}_{i+1})$, conditional on the value of
${\cal R}^p(B^{p,\vareps}_i)={\cal R}^p(b)$, is then
\[
{\cal R}^p\bigl(B^{p,\vareps}_{i+1}\bigr)= \cases{
{\cal R}^p(b)^*+N^{p,\vareps}\delta_{\gamma_p^{-1}A^{p,\vareps}}, &\quad if
$A^{p,\vareps}<a_{1}^*$ and $A^{p,\vareps}\neq
a_{1}$,
\cr
{\cal R}^p(b)^*, &\quad otherwise, }
\]
where we define ${\cal R}^p(a_1):=\supp({\cal R}^p(b))$, ${\cal
R}^p(b)^*:={\cal R}^p(b)-\delta_{{\cal R}^p(a_1)}$ and ${\cal
R}^p(a_1^*):=\supp({\cal R}^p(b)^*)$.

If we can now show convergence in distribution of
\[
\bigl(\gamma_p^{-1}A^{p,\vareps},N^{p,\vareps}
\bigr) \stackrel{d}{=} \bigl(\gamma_p^{-1}A_1,
\zeta'_{A_1} \bigr) \mid A_1\ge
\gamma_p\vareps
\mathop{\longrightarrow}_{p\to\infty}
\bigl(A_1^\vareps,N_1^\vareps\bigr),
\]
then our claim on point process convergence will follow from the
description of the transition rule for $(B^\vareps_i; i\ge1)$ from
Theorem~\ref{thmlawofcontscoal} and the standard convergence
arguments for a sequence of Markov chains based on weak convergence of
their initial values and transition laws.

We now express $(A_1,\zeta'_{A_1})$ in terms of the great-aunt measure
$\rho^{(p)}$ of a single quasi-stationary BGW tree with offspring
distribution $\xi_p$.
Consider the $(0,1)$ individual in our doubly infinite embedding of
quasi-stationary BGW genealogies. Recall (see Remark~\ref{rembinom})
that $(\zeta'_n; n\ge1)$ is a sequence of independent r.v.s which
conditionally on $(\rho^{(p)}_n, n\ge1)$ are binomial with parameters
$\rho^{(p)}_n$ and $p_{n-1}=\PP(Z^{(p)}_n\neq0\mid Z^{(p)}_0=1)$.

First for $A_1$, take any $x>\vareps$, then
\begin{eqnarray*}
\PP(A_1>\gamma_p x\mid A_1\ge
\gamma_p\vareps)&=&\prod_{i=\lceil\gamma
_p\vareps
\rceil}^{\lfloor\gamma_p x\rfloor}
\PP\bigl(\zeta'_i=0\bigr)=\EE\prod
_{i=\lceil
\gamma_p\vareps\rceil}^{\lfloor\gamma_p x\rfloor
}(1-p_{i-1})^{\rho
^{(p)}_i} \\
&=&\EE
\exp\Biggl(\sum_{i=\lceil\gamma_p\vareps\rceil
}^{\lfloor
\gamma_p x\rfloor}
\rho_i^{(p)}\ln(1-p_{i-1}) \Biggr).
\end{eqnarray*}
Let $\tau^{(p)}_{\mathrm{ext}}$ denote the extinction time of a BGW
process $Z^{(p)}$ started with $Z_0^{(p)}=p$, then
\[
1-p_{n-1}=\PP\bigl(Z^{(p)}_{n-1}=0\mid Z^{(p)}_0=1
\bigr)=\PP\bigl(\tau^{(p)}_{\mathrm{ext}}\le n-1\bigr)^{1/p.}
\]
%
Let $\tau_{\mathrm{ext}}$ denote the extinction time of the CSB process
$Z$ started with $Z_0=1$, then the assumption $\mathop{\lim\inf
}_{p\to
\infty}\PP(Z^{(p)}_{[\delta\gamma_p]}=0)>0$ 
guarantees that whenever \mbox{$\gamma_p^{-1}i_p\to i$}, we have
%
\[
(1-p_{_{i_p-1}})^p=\PP\bigl(\gamma_p^{-1}
\tau^{(p)}_{\mathrm{ext}}\le\gamma_p^{-1}i_p-
\gamma_p^{-1}\bigr)\mathop{\longrightarrow}_{p\to\infty}
\PP(\tau_{\mathrm{ext}}< i)=e^{-v(i)}.
\]
%
Let us define $f^{\vareps,x}(\cdot)= v(\cdot)\mathbf{1}_{[\vareps
,x]}(\cdot)$ and a sequence of functions $\{f_p\}$ such that
\[
f^{\vareps,x}_{p,\gamma_p}(\cdot)=f^{\vareps,x}_p\bigl(
\gamma_p^{-1}\cdot\bigr)= - \ln(1-p_{_{\cdot-1}})^p
{\bolds\chi}^{\vareps,x}_p(\cdot),
\]
where ${\bolds\chi}^{\vareps, x}_p$ is a sequence of bounded
continuous functions approximating $\mathbf{1}_{[\lceil\gamma
_p\vareps
\rceil,\lfloor\gamma_px\rfloor]}$ converging pointwise to $\mathbf
{1}_{[\vareps,x]}(\cdot)$.
Then we have that $f^{\vareps,x}_p(\gamma_p^{-1} i_p)\to f^{\vareps
,x}(i)= v(i)\mathbf{1}_{\vareps\le i\le x}$ as $p\to\infty$ whenever
$\gamma_p^{-1}i_p\to i$. By Theorem~\ref{thmpropertiesrho}(iii)
it follows that
\[
\sum_{i=\lceil\gamma_p\vareps\rceil}^{\lfloor\gamma_p x\rfloor} \frac
{\rho_i^{(p)}}{p}
\ln(1-p_{_{i-1}})^p \mathop{\longrightarrow
}^{d}_{p\to
\infty} -\bigl\langle\rho^0,v
\mathbf{1}_{[\vareps,x]}\bigr\rangle
\]
showing that
\[
\PP(A_1>\gamma_p x\mid A_1\ge
\gamma_p\vareps)\mathop{\longrightarrow}_{p\to
\infty} \EE\exp
\bigl(- \bigl\langle\rho^0,v\mathbf{1}_{[\vareps,x]}\bigr\rangle
\bigr).
\]
By Definition~\ref{dfnrho} and equation (\ref{eqAepsfromPP}) we have
\begin{eqnarray*}
&&
\EE\exp\bigl(-\bigl\langle\rho^0,v\mathbf{1}_{[\vareps,x]}\bigr
\rangle\bigr)\\
&&\qquad=\exp- \biggl(\beta\int_{[\vareps, x]}du\,v(u)+\int
_{[\vareps, x]}du\int_{[0,\infty)}\pi^{(u)}(dr)
\bigl(1-e^{-rv(u)} \bigr) \biggr),
\end{eqnarray*}
which together with the results of Theorem~\ref{thmlawofAeps1}
proves that $\gamma_p^{-1}A_1\mid A_1>\gamma_p\vareps
\stackrel{d}{\longrightarrow} A_1^\vareps$.

We next show that for any $y>\vareps$ the sequence $\sum_{i=\gamma
_p\vareps}^{\gamma_py}\zeta'_i$ converges as $p\to\infty$ to a Cox
point process on $[\vareps,y]$ whose intensity measure given $\rho^0$
is $\langle\rho^0,v\mathbf{1}_{[\vareps,y]}\rangle$.
For any $\lambda>0$,
\begin{eqnarray*}
\EE\exp\Biggl(-\lambda\sum_{i=\lceil\gamma_p\vareps\rceil
}^{\lfloor
\gamma_p y\rfloor}
\zeta'_i \Biggr)&=&\EE\prod_{i=\lceil\gamma
_p\vareps\rceil
}^{\lfloor\gamma_p y\rfloor}
\bigl(1-p_{i-1}+p_{i-1}e^{-\lambda}\bigr)^{\rho
^{(p)}_i}\\
&=&
\EE\exp\Biggl(\sum_{i=\lceil\gamma_p\vareps\rceil}^{\lfloor
\gamma_p
y\rfloor}
\rho^{(p)}_i \ln\bigl(1-p_{i-1}+p_{i-1}e^{-\lambda}
\bigr) \Biggr).
\end{eqnarray*}
Let us define $g^{\vareps,y}(\cdot)= v(\cdot)(1-e^{-\lambda
})\mathbf
{1}_{[\vareps,y]}(\cdot)$ and a sequence of functions $\{ g^{\vareps
,y}_p\}$ such that
\[
g^{\vareps,y}_{p,\gamma_p}(\cdot)=g^{\vareps,y}_p\bigl(
\gamma_p^{-1}\cdot\bigr)= -\ln\bigl(1-p_{_{\cdot-1}}+p_{_{\cdot
-1}}e^{-\lambda}
\bigr)^p {\bolds\chi}^{\vareps,y}_p(\cdot).
\]
%
Since whenever $\gamma_p^{-1}i_p\to i$, we have
\begin{eqnarray*}
p_{_{i_p}}&=&\PP\bigl(Z^{(p)}_{i_p}\neq0\mid Z^{(p)}_0=1
\bigr)=1-\PP\bigl(\tau_{\mathrm{ext}}^{(p)}\le i_p
\bigr)^{1/p}\\
&\approx&1-\PP\bigl(\tau_{\mathrm{ext}}< \gamma_p^{-1}i_p
\bigr)^{1/p}\approx v\bigl(\gamma_p^{-1}i_p
\bigr)/p,
\end{eqnarray*}
it follows that $g^{\vareps,y}_p(\gamma_p^{-1}i_p)\to g^{\vareps
,y}(i)=v(i)(1-e^{-\lambda})\mathbf{1}_{\vareps\le i\le y}$ as $p\to
\infty$ whenever $\gamma_p^{-1}i_p\to i$.
By Theorem~\ref{thmpropertiesrho}(iii),
\[
\EE\exp\Biggl(-\lambda p^{-1} \sum_{i=\lceil\gamma_p \vareps\rceil
}^{\lfloor\gamma_p y\rfloor}
\zeta'_i \Biggr) \mathop{\longrightarrow}_{p\to\infty}
\EE\exp\bigl( \bigl\langle\rho^0,v\mathbf{1}_{[\vareps,y]}\bigr\rangle
\bigl(e^{-\lambda}-1\bigr) \bigr),
\]
showing that (e.g.,~\cite{Kal}, Theorem 16.29)
%
\begin{equation}
\label{eqconvtoCox} \sum_{i=\lceil\gamma_p \vareps\rceil}^{\lfloor\gamma
_p y\rfloor
}
\zeta'_i \mathop{\longrightarrow}^{d}_{p\to\infty}
\Xi_{[\vareps,y]},
\end{equation}
where $\Xi$ is a Cox point process whose intensity measure
conditionally on $\rho^0$ is $\langle\rho^0,v\rangle$.
Finally, take $\vareps<x'<x<x''<y$, such that $\gamma_p(x'-x)\to0$ and
$\gamma_p(x''-x)\to0$
as $p\to\infty$.
\begin{eqnarray*}
&&\PP\bigl(\gamma_p x'\le A_1\le
\gamma_p x'',\zeta'_{A_1}=n\mid
A_1>\gamma_p\vareps\bigr)
\\
&&\qquad = \PP\bigl(\gamma_p x'\le A_1\le
\gamma_p x''\mid A_1>
\gamma_p\vareps\bigr) \PP\bigl(\zeta'_{A_1}=n\mid
\gamma_p x'\le A_1\le\gamma_p
x''\bigr)
\\
&&\qquad\mathop{\longrightarrow}_{p\to\infty}\PP\bigl(A_1^\vareps
\in dx\bigr) \PP\bigl(\Xi_{[\vareps,y]}(dx)=n \mid \Xi_{[\vareps
,y]}(dx)\neq
\varnothing\bigr).
\end{eqnarray*}
%
%
Evaluating the probability that the point process with intensity
measure $\langle\rho^0,v\rangle$ takes on values $n=1$ and $n\ge2$ at
height $x$,
gives precisely the formulae given in Theorem~\ref{thmlawofAeps1},
showing that for all $n\ge1$ and $x>\vareps$,
\[
\PP\bigl(A_1\in\gamma_p\,dx, \zeta'_{A_1}=n\mid A_1>
\gamma_p\vareps\bigr)\mathop{\longrightarrow}_{p\to\infty} \PP
\bigl(A_1^\vareps\in dx, N_1^\vareps=n
\bigr),
\]
%
%
and the proof is complete.
\end{pf}

\section{Two applications in the discrete case}

In this section, we come back to the discrete case to display two
further results on the coalescent point process.

\subsection{The linear-fractional case}\label{subseclinfrac}

A BGW process is called \textit{linear fractional} if there are two
probabilities $a$ and $b$ such that
\[
f(s)=a+\frac{(1-a)(1-b)s}{1-bs},\qquad s\in[0,1].
\]
In other words, $\xi$ is product of a Bernoulli r.v. with parameter
$(1-a)$ and a geometric r.v. with parameter $(1-b)$ conditional on
being nonzero. 
The expectation $m$ of $\xi$ is equal to
\[
m=\frac{1-a}{1-b},
\]
so that this BGW process is (sub)critical iff $a\stackrel{(>)}{=}b$.
In general, the coalescent point process $(A_i,i\ge1)$ is not itself a
Markov process, but in the linear-fractional case, it is a sequence of
i.i.d. random variables. An alternative formulation of this observation
was previously derived in~\cite{Ran}.
%
\begin{prop}
\label{propcpplf} In the linear-fractional case with parameters
$(a,b)$, the branch lengths of the coalescent $(A_i,i\ge1)$ are i.i.d.
with distribution given by
\[
\PP(A_1> n)= \frac{b-a}{bm^n-a},\vadjust{\goodbreak}
\]
when $a\not=b$, and when $a=b$ (critical case), by
\[
\PP(A_1> n)= \frac{1-a}{na+1-a}.
\]
\end{prop}
\begin{pf}
Recall from Theorem~\ref{thmcpp} that $A_i=\min
\{
n\ge1\dvtx  D_i(n)\not=0\}$, so in particular, we can set $A_0:=+\infty$.
We prove by induction on $i\ge1$ the following statement $({\cal
S}_i)$. The random variables $(D_i(n);n\ge1)$ are independent r.v.s
distributed as $\zeta_n'$, and independent of $(A_0,\ldots, A_{i-1})$.
Observe that $({\cal S}_1)$ holds thanks to Theorem~\ref{thmcpp}.
Next, we let $i$ be any positive integer, we assume $({\cal S}_i)$ and
we prove that $({\cal S}_{i+1})$ holds. Elementary calculus shows that
$\zeta_n$ has a linear-fractional distribution, so that $\zeta_n'$ has
a geometric distribution. We are now reasoning conditionally given
$A_i=h$. By $({\cal S}_i)$ and the definition of $A_i$, we get that
conditional on $A_i=h$:
\begin{itemize}
\item the r.v.s $(D_i(n);n>h)$ are independent r.v.s distributed as
$\zeta_n'$, and independent of $(A_0,\ldots, A_{i-1})$,
\item
$D_i(h)$ has the law of $\zeta_h'$ conditional on $\zeta_h'\not=0$, and
it is independent of $(D_i(n);n>h)$ and $(A_0,\ldots, A_{i-1})$,
\item
$D_i(n)=0$ for all $n<h$.
\end{itemize}
Let us apply the transition probability defined in the theorem.
First, $D_{i+1}(n)=D_i(n)$ for all $n>h$, so the r.v.s
$(D_{i+1}(n);n>h)$ are independent r.v.s distributed as $\zeta_n'$,
and independent of $(A_0,\ldots, A_{i-1})$. Second,
$D_{i+1}(h)=D_i(h)-1$ has the law of $\zeta_h'-1$ conditional on
$\zeta_h'\not=0$, which is the law of $\zeta_h'$ because $\zeta_h'$ is
geometrically distributed, and it is independent of $(D_i(n);n>h)$ and
$(A_0,\ldots, A_{i-1})$. Third, the r.v.s $(D_{i+1}(n);n<h)$ are new
independent r.v.s distributed as $\zeta_n'$. As a consequence,
conditional on $A_i=h$, the r.v.s $(D_{i+1}(n);n\ge1)$ are
independent r.v.s distributed as $\zeta_n'$, and independent of $(A_0,\ldots, A_{i-1})$. Integrating over $h$ yields $({\cal S}_{i+1})$.

We deduce from $({\cal S}_{i})$ that $A_i$ is independent of $(A_0,\ldots, A_{i-1})$ and is distributed as $A_1$. The computation of the
law of $A_1$ stems from well-known formulae involving linear-fractional
BGW processes (see~\cite{AN}), namely,
\[
f_n(s) = 1-\frac{(1-(a/b))(1-s)}{m^{-n}(s-(a/b))+1-s},
\]
if $a\not=b$, and
\[
f_n(s) = \frac{na-(na+a-1)s}{1-a+na-nas},
\]
when $a=b$. Indeed, it is then straightforward to compute
\[
f_n'(0) = m^{-n}\frac{(b-a)^2}{(b-am^{-n})^2}
\quad\mbox{and}\quad 1-f_n(0) = \frac{b-a}{b-am^{-n}},
\]
if $a\not=b$, whereas
\[
f_n'(0) = \frac{(1-a)^2}{(1-a+na)^2} \quad\mbox{and}\quad
1-f_n(0) = \frac
{1-a}{1-a+na},
\]
when $a=b$. Thanks to Theorem~\ref{thmcpp}, the ratio of these
quantities is \mbox{$\PP(A_1> n)$}.
\end{pf}
%
%

\subsection{Disintegration of the quasi-stationary distribution}
\label{subsecdisintegration}
In this subsection, we assume that $f'(1)<1$ (subcritical case). Then
it is well known~\cite{AN} that there is a probability $(\alpha
_k)_{k\geq1}$ with generating function, say $a$,
\[
a(s) = \sum_{k\geq1} \alpha_k
s^k,\qquad s\in[0,1],
\]
such that
\[
\lim_{n\tendinfty}\PP({Z_n}=k\mid Z_n \ge1) =
\alpha_k,\qquad k\ge1.
\]
This distribution is known as the \textit{Yaglom limit}. It is a
quasi-stationary distribution, in the sense that
\[
\sum_{k\geq1} \alpha_k
\PP_k(Z_1 = j\mid Z_1\not=0)=
\alpha_j,\qquad j\geq1.
\]
%

Set
\[
U:=\min\{i\ge1\dvtx  A_i =+\infty\}.
\]
Then for any $i\le U <j$, the coalescence time between individuals
$(0,i)$ and $(0,j)$ is $\max\{A_k\dvtx i\le k<j\}=+\infty$, so that $(0,i)$
and $(0,j)$ do not have a common ancestor. Now set
\[
V:=\max\{A_k\dvtx  1\le k <U\},
\]
the coalescence time of the subpopulation $\{(0,i)\dvtx 1\le i \le U\}$
(where it is understood that $\max\varnothing=0$), that is, $-V$ is
the generation of the most recent common ancestor of this subpopulation.
We provide the joint law of $(U,V)$ in the next proposition.
%
\begin{prop}
\label{propmrcaqsd}
The law of $V$ is given by
\[
\PP(V\le n) = \alpha_1 \frac{1-f_n(0)}{f_n'(0)}= \frac{\alpha_1}{\PP
(Z_n=1\mid Z_n\not=0)},\qquad n
\ge0.
\]
Conditional on $V=n\ge1$, $U$ has the law of $Z_n$ conditional on
$\zeta_n\ge2$.
In addition, $U$ follows Yaglom's quasi-stationary distribution, which
entails the following disintegration formula:
\[
a(s)=\PP(V=0)s+\sum_{n\ge1}\PP(V=n) \EE
\bigl(s^{Z_n}\mid\zeta_n \ge2\bigr),\qquad s\in[0,1].
\]
\end{prop}
%
\begin{rem}
In the linear-fractional case with parameters $(a,b)$ $(a>b$,
subcritical case), the Yaglom quasi-stationary distribution is known
to be geometric with failure probability $b/a$. Thanks to Proposition
\ref{propcpplf}, the branch lengths are i.i.d. and are infinite
with probability $1-(b/a)$. Then thanks to the previous proposition,
the quasi-stationary size $U$ is the first $i$ such that $A_i$ is
infinite, which indeed follows the aforementioned geometric law.
\end{rem}
\begin{pf*}{Proof of Proposition~\ref{propmrcaqsd}}
From the transition probabilities of the Markov
chain $(D_i;i\ge1)$ given in Theorem~\ref{thmcpp}, we deduce that
\[
V=\max\bigl\{n\dvtx D_1(n)\not=0\bigr\},
\]
so that, thanks to Theorem~\ref{thmcpp},
\[
\PP(V\le n)=\prod_{k\ge n+1} \PP\bigl(\zeta_k'=0
\bigr)=\frac{\prod_{k\ge1} \PP
(\zeta_k'=0)}{\prod_{k=1}^n \PP(\zeta_k'=0)}=\frac{\PP(A_1=+\infty)}{\PP(A_1>
n)}.
\]
Now, since
\[
\PP(A_1> n) = \frac{f_n'(0)}{1-f_{n}(0)}= \PP(Z_n=1\mid
Z_n\not=0),
\]
and because this last quantity converges to $\PP(A_1=+\infty)=\alpha_1$
as $n\tendinfty$, we get the result for the law of $V$.

Recall that $(-n,\mathfrak{a}_1(n))$ is the ancestor, at generation
$-n$, of $(0,1)$, so that the total number of descendants $\Upsilon
_n:=Z^{(-n,\mathfrak{a}_1(n))}(n)$ of $(-n,\mathfrak{a}_1(n))$ at
generation $0$ has the law of $Z_n$ conditional on $Z_n\not=0$. Now,
since no individual $(0,j)$ with $j>U$ has a common ancestor with
$(0,1)$, we have the inequality $\Upsilon_n \le U$. On the other hand,
$(0,U)$ and $(0,1)$ do have a common ancestor [all coalescence times
$(A_i;1\le i\le U-1)$ are finite], so that there is $n$ such that
$\Upsilon_n=U$. Since the sequence $(\Upsilon_n)$ is obviously
nondecreasing, it is stationary at $U$, that is,
\[
\lim_{n\tendinfty} \Upsilon_n = U \qquad\mbox{a.s.}
\]
Since $\Upsilon_n$ has the law of $Z_n$ conditional on $Z_n\not=0$, $U$
follows the Yaglom quasi-stationary distribution.

Actually, since $-V$ is the generation of the most recent common
ancestor of $(0,1)$ and $(0,U)$, we have the following equality:
\[
U=\Upsilon_{V}.
\]
Now recall that $V=\max\{n\dvtx D_1(n)\not=0\}$. By definition of $D_1$, we
can write $\{V=n\}= E_n\cap F_n$, where $E_n:=\{D_1(k)=0, \forall k>n\}
$ and
\[
F_n:=\bigl\{D_1(n)\ge1\bigr\}=\bigl\{\# {\cal D}(n,1)
\ge2\bigr\}.
\]
Now observe that $E_n$ is independent of all events of the form
$F_n\cap\{\Upsilon_n=k\}$ (it concerns the future of parallel
branches), so that $\PP(\Upsilon_n=k\mid V=n)=\PP(\Upsilon_n=k\mid
F_n)$. In other words, $\Upsilon_n$ conditional\vadjust{\goodbreak} on $V=n$ has the law of
$Z_n'$ conditional on $\zeta_n\ge2$, where $Z_n'$ is $Z_n$ conditional
on $Z_n\not=0$. Since $\{Z_n\not=0\}=\{\zeta_n \ge1\}$, we finally get
that conditional on $V=n$, $U$ has the law of $Z_n$ conditional on
$\zeta_n\ge2$.
%
%
\end{pf*}

\section*{Acknowledgments}

A. Lambert wishes to thank Julien Berestycki and Olivier H\'{e}nard for
some interesting discussions on the topic of this paper. A. Lambert and
L. Popovic thank the referees for suggestions that improved the
exposition of the paper.



\printaddresses

\end{document}